\documentclass[10pt,leqno] {amsart}
\makeatletter
\@namedef{subjclassname@2020}{2020 AMS  Subject Classification}
\makeatother
 \usepackage{mathrsfs,amssymb, color}

\usepackage{enumitem}
 \usepackage{graphicx}

\newtheorem{thm}{Theorem} [section]
\newtheorem{thm*}{Theorem*} [section]
\newtheorem{lem}[thm]{Lemma}

\newtheorem{corol}[thm]{Corollary}
\newtheorem{corol*}[thm]{Corollary*}
\newtheorem{prop}[thm]{Proposition}
\newtheorem{defin}[thm]{Definition}

\newtheorem{rem}[thm]{Remark}
\newtheorem{rems}[thm]{Remarks}

\newtheorem{ssect}[thm]{}

\renewcommand\P{\mathbb P}
\renewcommand\AA{{\mathcal A}}
\newcommand \BB{{\mathcal B}}

 \newcommand \DD{{\mathcal D}}

\newcommand \FF{{\mathcal F}}
\newcommand \GG{{\mathcal G}}
\newcommand \HH{{\mathcal H}}

\newcommand\JJ{{\mathcal J}}
 \newcommand \KK{{\mathcal K}}
\newcommand \PP{{\mathcal P}}

\newcommand\TT{{\mathcal T}}
 \newcommand \UU{{\mathcal U}}

 \newcommand \I{{\mathbb I}}

\newcommand\QQ{{\mathbb Q}}
\newcommand\R{{\mathbb R}}
  
 \newcommand\ZZ{\mathscr  Z}
  \newcommand\M{\mathcal  M}
  \newcommand\T{\mathscr T}
 \newcommand\X{\mathscr X}

\newcommand\A{\boldsymbol A}
\newcommand\B{\boldsymbol B}
\newcommand\C{\boldsymbol C}
\newcommand\D{\boldsymbol D}
\newcommand \Uu{\boldsymbol U}
\renewcommand \ss{\boldsymbol \s}
 \def\SSS{\boldsymbol\Sigma}
 \renewcommand \SS{\boldsymbol S}
\renewcommand \R{\boldsymbol R}
\newcommand \Q{\boldsymbol Q}

 \newcommand \borel{\boldsymbol{\varDelta}^1_1}
 
 \newcommand \borm[1]{\boldsymbol{\Pi}^0_{#1}}
 \newcommand\ana{\boldsymbol\Sigma^1_1}
 \newcommand\ca{\boldsymbol\Pi^1_1}
\newcommand \Gd{\boldsymbol{G}_\delta}

\newcommand\cospi{\check\AA (\ca)}
 \newcommand\spi{\AA (\ca)}

\newcommand\cpca{\boldsymbol\Pi^1_2}
 
\newcommand\pca{\boldsymbol\Sigma^1_2}

\newcommand\baire{ \omega^\omega}
\newcommand\G{\boldsymbol \Gamma}

\let\a=\alpha

\let\wo=\omega
\let\s=\sigma
\let\d=\delta
\let\e=\varepsilon
\let\eps=\e

\let\th=\theta
\let\Th=\Theta
\def\f{\mathfrak f}

 \let\ph=\varphi

\renewcommand\S{\Sigma}
\newcommand\seq{\omega^{<\omega}}
   
\def\i{\ifmmode \mathchar "231\else \char "10\fi}
\def\rond{{\scriptstyle\circ}}
 
\DeclareMathOperator\diam{\rm diam}

\DeclareMathOperator\dom{\rm dom}
\DeclareMathOperator\card{\rm card}

\let\<=\langle  \let\>=\rangle
\let\nd=\noindent
\let\ct=\centerline
\let\mk=\medskip
\let\sk=\smallskip

\def\sn{\nobreak\smallskip\par}

\let\mns=\setminus

\newcommand \Ind[2]{  {#1}^{\{#2\}}  }

\newcommand\cat{{}^\frown}      

\newcommand\td[1]{\tilde #1}
\newcommand\Ex{\exists \,}
\newcommand\Vx{\forall \,}
\newcommand\m{^{-1}}

 \let\0=\emptyset
\newcommand\nv{\neq\emptyset}
 
\newcommand\cc{\overset  \circ}
 \newcommand\cp{\;\hat  }

 \newcommand\si{\ {\rm if} \ }

 \newcommand\et{\ {\rm and} \ }
 \newcommand\ett{\quad {\rm and} \quad}
 
  \newcommand\ou{\ {\rm or} \ }

 \newcommand\adh[1]{\overline{#1}}

 \def\precneq{ \mathrel{ \mathop{\kern0pt\smash\prec}\mkern -16mu_{\neq}}}
 
\def\abs#1{\left\vert#1\right\vert}

\def\hom@{homeomorphism}

 \let\<=\langle
\let\>=\rangle

\def\st{such that }

\def\Carc{\mathcal C_{\rm arc}} 
 
\def\Cloc{\mathcal C_{\rm loc}}

 \newcommand\pcon[1]{{\rm PCON}_{#1}}
  \newcommand\dg{degenerated}

  \newcommand\SaF{S_{(\a,F)}}

\def\pcon{\Carc(\R^2)}

 \subjclass[2020]{54H05,  03E15}

 \title[Compact sets in the plane]{The descriptive complexity  of the set of  \\ arc-connected compact subsets of the plane}

\author{Gabriel Debs} 
\author{Jean Saint Raymond}

\address{Gabriel Debs, Sorbonne Universit\'e, Universit\'e Paris Diderot, CNRS, Institut de Math\'ematiques de Jussieu-Paris Rive Gauche, 
IMJ-PRG, F-75005, Paris, France}
\email{gabriel.debs@imj-prg.fr}

\address{Jean Saint Raymond, Sorbonne Universit\'e, Universit\'e Paris Diderot, CNRS, Institut de Math\'ematiques de Jussieu-Paris Rive Gauche, IMJ-PRG, F-75005, Paris, France}
\email{jean.saint-raymond@imj-prg.fr}

\begin{document}

\begin{abstract} We compute the exact descriptive class of the  set of all compact arc-connected  subsets of $\R^2$, which turns out to be strictly higher than  the classical $\ana$ and $\ca$ classes  of analytic and coanalytic sets, but strictly lower than  the  class $\boldsymbol \Pi^1_2$ which is the exact descriptive class   of the  set of  all compact  arc-connected subsets of $\R^3$.
\end{abstract} 

 \maketitle

If $X$ is any Polish  space then it
 follows readily from the  definitions that the set  
$\KK\Carc(X)$ 
 of all compact arc-connected  subsets of $X$,
 is a  $\boldsymbol \Pi^1_2$ subset of the   space   $\KK(X)$ 
  of all compact subsets of $X$, endowed with the Vietoris topology. Moreover   Ajtai and Becker  showed   independently (see \cite{kc},   Theorem 37.11)  that  
   the  set  $\KK\Carc(\R^3)$  
    is actually $\boldsymbol \Pi^1_2$-complete.

\sk

The  goal of the present work is to compute the exact descriptive complexity of the  set  $\KK\Carc(\R^2)$.
More generally, given any space $X$ we consider the set $\Carc(X)$ of  all arc-connected closed  subsets of $X$, viewed as a subset of the  space  $\FF(X)$  of all closed subsets of $X$, endowed with the 
Effros Borel structure.  

By a {\em planar Polish space} we mean a subspace of $\R^2$, on which the induced  topology is Polish, that is  a $\Gd$ subset of $\R^2$.
Our    first main result is the following:

\sk

\nd{\bf Theorem A.} {\it  For any planar Polish space $X$ the set $\Carc(X)$ is a $\cospi$ set.}
 
\mk

\nd  where $ \cospi$  denotes the class of all complements of sets obtained from $\ca$ sets by  Suslin operation $\AA$. Let us recall  that it was already known from Ajtai and Becker  work that the  set $\KK\Carc(\R^2)$ is not $\ca$ (see \cite{kc}).

The proof of  Theorem A   relies on recent results from \cite{dsr}. Also as a by-product  of this proof we obtain  the following   property  which is  specific  of the plane topology, since the analog is no more true in $\R^3$.

\sk

\nd{\bf Theorem B.} {\it  For any planar arc-connected Polish space $X$ there exists a        Borel mapping
 which to any pair $(x,y)$ of distinct points in $X$ assigns an arc $J\subset X$  with endpoints $x$ and $y$.}

\sk 
   
In fact  the proof of Theorem A relies on a parametrized version of Theorem B in which the space $X$ is replaced by a variable closed subset of $X$, given with some code in an auxiliary space. The precise statement of this latter result  (Theorem \ref{main})  necessitates a number of preliminaries,  and we refer the reader to Section \ref{unifcode} for more details.

\sk 

The second main result is that   if  $X=\R^2$,  or the unit square $\I^2$, then the complexity bound   given by  Theorem~A  is  best possible.   More precisely we prove:
\sk

 \nd{\bf Theorem C.} {\it  The set  $\Carc(\I^2)$ is $\cospi$-complete.}
 
 \sk
 
 We also give in Section \ref{unifcode} several cases in which the set $\Carc(X)$ is in a strictly smaller class than the class $\cpca$.  
 
 It is worth noting  that the  class $\cospi$  appeared already in previous complexity computations in the hyperspace $\FF(X)$. We mention   the following two results  from  \cite{dsr0} (Theorems 6.3 and 6.4) where $\Cloc(X)$ denotes the set of all locally connected closed subsets of $X$:   

 \sk

{\it  a) For any  Polish space $X$ the set $\Cloc(X)$ is a $\cospi$ set.

b) There exists a    Polish space $X\subset \I^3$ for which  the set $\Cloc(X)$ is $\cospi$-complete}

 \sk
 
\nd  Note however that unlike for $\Carc(X)$, if $X$ is a compact space   the set $\Cloc(X)$ is Borel (\cite{dsr0} Proposition 6.1).

\section{Some descriptive preliminaries}

Throughout this work by a  {\it ``space"}   we shall always mean a subset of some Polish space,
though we shall introduce in some situations   an additional   (possibly non separable) metric or topology on the initial given space. However all descriptive notions  we shall consider   will  always refer   to the Borel structure inherited from the Polish topology.

  We shall consider  various classical, and less classical, descriptive classes. 
 For such a class $\G$ we denote by $\check \G$ its  dual class,  that is the class  of all complements of sets in 
 $\G$, and by $\Delta(\G)= \Delta(\check\G)$ the self-dual class  $\G\cap\check\G$.
Then  given any   space $X$,   subspace of some Polish space $\td  X$, we denote by $\G(X)$ the set of all subsets of $X$ which are the trace on $X$ of a subset of $\td X$ which is in $\G$. If the class $\G$ is  closed under Borel isomorphim, 
which will always be the case in this work, then  $\G(X)$ does not depend on the particular choice of the surrounding space $\td X$.  Note   that  in general $\G(X)$ is a proper subset of $\PP(X)\cap \G(\td X)$, and
  $\Delta(\G)(X)$ is a proper subset of  $\Delta(\G(X))$.

For any class $\G$ we   consider also the
  class $\AA(\G)$   obtained from $\G$ by    Suslin operation $\AA$  and  denote by
$\check \AA(\G)$ its  dual class.  Since    countable unions   and countable    intersections  are particular instances of    operation $\AA$,  it  follows   from the idempotence of operation  $\AA$ that the classes $ \AA(\G)$ and $\check \AA(\G)$ are closed under countable unions and intersections.

  Following  logician notation we denote by  $\boldsymbol\Sigma^1_n$ and $\boldsymbol\Pi^1_n$ the classical projective classes, and set  $\boldsymbol\Delta^1_n=\Delta(\boldsymbol\Sigma^1_n)=\Delta(\boldsymbol\Pi^1_n)$. In particular   $\ana$ and $ \ca$  denote  respectively   the classes of analytic and  coanalytic  sets,   $\pca$ the class of projections of $\ca$  sets  and   $\cpca$  its dual class.

\sk

We  recall that a mapping $f:X\to Y$     is said to be {\em $\G$-measurable} if the inverse image of any open subset of $Y$ is in $\G(X)$. If moreover $\G$ is closed under countable unions and intersections then the inverse image of any Borel subset of $Y$   is in $\G$.     It follows that  for any class $\G$, the notions of  $\AA(\G)$-measurability and $\check \AA(\G)$-measurability are the same. 
Also since $\ana=\AA(\borel)$ and $\ca=\check\AA(\borel)$ we have   
the following property:

\begin{prop} \label{A(G)}
 The inverse image of any $\ana$ (respectively $\ca$) set  by an $\AA(\G)$-measurable mapping      is in $\AA(\G)$ (respectively $\check \AA(\G)$). 
\end{prop}  

  In particular the right composition $f\rond g$   of an $\AA(\G)$-measurable mapping $f$ with a Borel mapping $g$, is $\AA(\G)$-measurable; and if $\G$ is closed under Borel isomorphisms then  the left  composition $h\rond f$   of  $f$ with a Borel mapping $h$ is $\AA(\G)$-measurable too.
  
Of particular interest for our study is  the   notion of bianalyticity that we recall.
  
\begin{defin} 
Given spaces $X,Y$:

  A subset $A$ of $X$ is said to be {\em bianalytic} in $X$ if   $A$ is in  
  $\ana(X)\cap \ca(X)$.

  A mapping $f:X\to Y$  is said to be {\em bianalytic} if  $f$ is $\ana$ (equivalently  $\ca$)-measurable.
\end{defin}

So the set of bianalytic subsets of $X$ is $\Delta(\ana(X))=\Delta(\ca(X))$ while the set of Borel subsets of $X$ is $\Delta(\ana)(X) =\Delta(\ca)(X)$.
Note that by  the separation Theorem of analytic sets,  if $X$ is analytic then the notions of 
bianalytic    and   Borel coincide. In fact  the notion of  bianalyticity is interesting mainly in the frame of  non analytic,  and more specifically coanalytic, spaces.  
Also in this latter context many properties of Borel sets (mappings) extend to  bianalytic sets (mappings), and we state next two such properties which we will need.

\begin{prop} \label{ext}  Let $X$ and $Y$ be coanalytic spaces. Any partial bianalytic mapping  $f:D\to Y$  with $D\subset X$,    admits an extension $\td f:\td D\to Y$ to a bianalytic mapping  with  coanalytic domain  $D\subset \td D \subset X$.
\end{prop}

The proof of Proposition \ref{ext} follows standard arguments.
Note that if moreover the graph of $f$ is a subset of some coanalytic set $Z\subset X\times Y$ then applying the previous result to the mapping $(Id_D, f): D\to Z$ we can impose that the graph of the mapping $\td f$ is also a subset of $Z$.

\sk

The second property we will need is less elementary and is a bianalytic replica of the classical fact that the projection of a Borel set with compact sections is Borel (a particular case of Arsenin-Kunugui Theorem).
The result is probably not new, but in the absence of known references we sketch the basic arguments of the proof.

\begin{thm}\label{unicity} For    any   space $X$, and any     standard Borel space $Y$, if    $B$  is  a bianalytic subset  of $ X\times Y$, with compact sections in $Y$ then the projection $P$ of $B$ on $X$ is  bianalytic,  the mapping   $\Phi: x\mapsto B(x)$ from  $P$  to $\KK(Y)$   is     bianalytic and admits a  bianalytic selection ({\it i.e.} there exists a   bianalytic mapping   $\ph: P\to Y$ \st for all $x\in P$, $\ph(x)\in B(x)$). 
\end{thm}

\begin{proof}  We embed $X$   into a  Polish space   $\td X$, and 
  fix  in $\td X\times  Y$ a  $\ana$  subset  $A$ and a  $\ca$ subset  $C$   \st  $
  B=A \cap(X\times Y) =C \cap(X\times Y)$. We also fix  on $Y$ a bounded distance $d$ and a basis 
  $(U_s)_{s\in S}$  of the topology, where $S\subset 2^{<\wo}$ is a tree and for all $s\in S$, 
  $\diam(U_s)< 2^{-\vert s\vert}$ and $\adh{ U_{s}} \subset U_{s^*}$ if $s\nv$  (where $s^*$ is the restriction of $s$ to $\vert s\vert -1 $); and we set $S_n=S\cap 2^n$.

  Let $\pi:\td X\times Y \to \td X$ be the canonical projection.
  Then for  all $s\in S$   the  set $P_s=\pi(A\cap (\td X\times U_s))$ is $\ana$, hence  
  $A_n=\bigcup_{s\in S_n}     P_s\times U_s $   is $\ana$ too, and $\td A=\bigcap_{n\in \wo}\downarrow A_ n$ is $\ana$ too. Moreover
  $$(x,y)\in  \td A \iff \Vx n, \  \Ex s\in S_n, \ y_n\in U_s, \ (x,y_n)\in A$$
and  since each section $B(x)$ of $B$ is compact,  then  $B\subset \td A$.

Then by the reduction property of the class $\ca$   we can fix a sequence $(C_n)_{n\in \wo}$ of pairwise disjoint  $\ca$ sets \st for all $n$, $C_n\subset (\td X\times Y) \mns A_n$ and  $\td C:=\bigcup_{n\in \wo}C_n =\bigcup_{n\in \wo} (\td X\times Y) \mns A_n = (\td X\times Y) \mns  \td A$.
So  $C'_n:=\bigcup_{m\leq n}C_m\subset (\td X\times Y) \mns A_n$  and $C''_n:=\td C\mns C'_n=\bigcup_{m>n}C_m$, hence $C'_n$ is a bianalytic subset of $\td C$, and $A_n\subset \td C_n= C''_n \cup \td C$.

Note that by definition  for all $s\in 2^n$:
$$P_s\subset \{x\in \td X: \ \{x\}\times U_s\subset A_n\}\subset \{x\in \td X: \ \{x\}\times U_s\subset \td C_n\}:=R_s$$ 
and $R_s$ is $\ca$. Then again by  the reduction property of the class $\ca$  we can find a bianalytic
set $P_s\subset Q_s\subset R_s$ \st $R_s\times U_s\subset \td C_n$.  

 Hence 
$B\subset \bigcap_n\bigcup_{s\in S_n} Q_s\times U_s$ and $P=\pi(B)\subset  X\cap \bigcap_n \bigcup_{s\in S_n} Q_s$.   Conversely if $x\in   X\cap \bigcap_n \bigcup_{s\in S_n} Q_s$ then it follows from  the compactness of $B(x)$   that $x\in P$, which proves that $P$ is a bianaytic  subset of $X$.
 Then for any open set $V$ in $Y$, since $B\cap (X \times V)$ is bianalytic in $X \times Y$,  the set
  $\pi(B\cap (X \times V)$ is a bianaytic  subset of $X$. It follows that the mapping $\Phi: P\to \KK(Y)$ is 
 bianaytic, and if ${\bf c}: \KK(Y)\to Y$ is any Borel choice mapping on $\KK(Y)$ ({\it i.e.} ${\bf c}(K)\in K$ for all $K\in \KK(Y)$) then
 $\ph={\bf c}\rond \Phi$ satisfies the conclusion of Theorem \ref{unicity}.
\end{proof}

\begin{ssect} The class  $\SSS$:
{\rm We denote by $\SSS$    the    smallest   class  containing both  classes  $\ana$  and $\ca$, and  closed under countable unions and intersections;
so $\SSS\subset \spi\cap\cospi$. 
We recall  the following classical result: 
  
\begin{thm} ({\sc Yankov - von Neumann}) \label{YvN} 
For any  $\ana$ set  $A\subset X\times Y$ in a product space, if $P$ is the projection of $A$ on $X$ then   there exists a  $\SSS$-measurable  
mapping $f: P\to Y$  with   graph  contained in $ A$.
 \end{thm} 
 } \end{ssect}
 
 We finish with  the following two  complements to Theorem \ref{YvN}.

 \begin{lem}
 Let $X$ and $Y$ be two Polish spaces and $\hat Y:=  Y\cup\{*\}$, where $ *\notin  Y$.
 Let  $D$ be a countable set and $f$ be a $\SSS$-measurable function from $X$ to  $\hat Y^D$.
If  for each $x\in X$, the set $\{ d\in D : f_d(x)\neq*\}$ is infinite, then there exists a $\SSS$-measurable function $\tilde f :  X\to  Y^\omega$ such that for all $x\in  X$
 $$\{ f_d(x) : d\in D\et f_d(x)\neq*\} = \{ \tilde f_n(x) : n\in \omega\}$$
\end{lem}
\begin{proof}
 Enumerate $ D$ as a sequence $(d_k)_{k\in \omega}$ and for all $x\in X$ denote $ D_x $ the countable  infinite set $ \{ d\in D : f_d(x)\neq *\}=\{ d\in D : f_d(x)\in  Y\}$.
 Then take for $ \tilde f_n(x)$ the $ n^{\rm th}$ term of the infinite sequence (with partial domain) $\bigl( f_{d_k}(x)\bigr)_{k\in D_x}$. This ensures that 
 $$ Y \supset
 \{ f_d(x) : d\in D\et f_d(x)\neq*\} =\{ f_d(x) : d\in D_x\} = \{ \tilde f_n(x) : n\in \omega\}$$
 
 To prove that $\tilde f$ is $\SSS$-measurable we have to show that for all $n\in \omega$ and all open set $ V\subset  Y$ the set $\tilde f_n\m(V)$ belongs to $ \SSS$. 
 For all $u\in \omega$ the set $  L_u=\{ x\in  X : f_{d_u}(x)\neq *\} $  belongs to $ \SSS$ and so does
  $ L_u^V=\{ x\in L_k : f_{d_u}(x)\in V\}$.
 Then  the set $ \{x : \tilde f_n(x)\in V\} $ is equal to $$
  \bigcup_{ u_1< u_2<\dots <u_n=k}\nolimits\Bigl( L_k^V\cap \bigcap_{i\leq n}\nolimits L_{u_i} \cap \bigcap _{i\in k\mns\{ u_1, u_2, \dots, u_n\}} \nolimits(X \mns L_{u_i})\Bigr)$$
 which belongs to $ \SSS$.
\end{proof}

\begin{lem} \label{compress}
 Let $ X$ and $ Y$ be Polish spaces, $*\notin Y$ and $ Z$ be an analytic subset of $  X\times  Y$. Then there exists a $ \SSS$-measurable fonction $ f =  X \to \hat Y=  Y\cup \{*\}$ \st $(x,f(x))\in Z $
 if $x$ belongs to the projection $\pi(Z) $ of $ Z$ on $ \X$ and
 $f(x)=* \iff x\notin \pi(Z)$.
\end{lem}
\begin{proof}
 It follows from Theorem \ref{YvN} that there is a $\SSS$-measurable fonction $ f :\pi(Z) \to  Y$ \st $f(x)\in \{ y\in  Y : (x,y)\in Z\}$ whenever $x$ belongs to the analytic set $ \pi(Z)$. 
 Then extending $ f$ by $ *$ on the coanalytic set $  X\mns \pi(Z)$ yields a $ \SSS$-measurable function from $ X$ to  $\hat Y$. 
 \end{proof}

 \section{The arc-connection relation} \label{s_ker}

In this section we present  briefly the main notions and known results,  which we will use freely in the sequel, concerning the arc-connection relation, namely in the plane. For more details we refer the reader to  \cite{dsr}.

  \begin{ssect}    Arcs: {\rm 
By an {arc} $I$ we mean  as  usual a compact space homeomorphic to the  unit interval $\I=[0,1]$. For an arc $I$  we    denote by $e(I)$ the set of its endpoints and set $\cc {I\;}=I\mns e(I)$.     
The set  $\JJ(X)$   of all arcs in some  space $X$, viewed as a subset  of  the space $\KK(X)$,    is Borel (in fact $\borm3=\boldsymbol F_{\s\d}$) and the mapping  $e:\JJ(X)$ to $\KK(X)$ is Borel (see \cite{bp}).

\sk 

  For practical reasons, we shall consider a   singleton $\{a\}$ as a {\em \dg\,arc}  with $a$ as a unique endpoint. All notions and notations for  arcs    extend trivially to \dg\, arcs.  
A set  will  said to be a {\em possibly \dg\,arc},   and we shall write {\em p.d.\,arc}, if it is either an   arc or a singleton. We denote by $\hat \JJ(X)$ the set   of all  p.d.\,arcs.
\sk

If $I$ is an arc, for any     $\{a,b\} \subset I$    we denote   by  $\Ind I {a,b}$ the sub-arc  of  $I$ with endpoints  $\{a,b\}$.
 The mapping    which to any arc $I$ and any  pair $\{a,b\}\subset I $ assigns the arc  $ \Ind I {a,b}$, is Borel,
 since its graph:
$\{(I,J)\in \JJ(X)^2: \  J\subset I \et  e(J)=\{a,b\}\}$ 
  is    Borel. 
}\end{ssect}

\begin{ssect} Triods:
{\rm 
 A {\em  simple triod}  in a  space $X$  is  a compact subset $T=J_0\cup J_1\cup J_2$ of $X$ which is the union of three  arcs $J_i$  such that for all  $ i\not=j$,  $J_i\cap J_j=\{c\}$ is a singleton. 
The arcs $J_i$, which  are  uniquely determined  up to a permutation,  are called the {\it branches} of $T$,    and $c$  is  called  the {\it center} of $T$. 

\sk

In this work we shall never consider the more general notion of {\it triod}   introduced initially by
Moore in \cite{mo}, and    in the sequel  by a  {\em  ``triod"}  we shall always mean a  {\em  ``simple triod".}

  \sk
  
Since the set $\JJ(X)$ of all arcs is a Borel subset of $\KK(X)$ and  the $\cup$ and $\cap$ operations on $\KK(X)$  are Borel,  it follows from the almost uniqueness of the decomposition of a simple triod, that if   $X$ is a  Polish space  then the set  $\TT(X)$   of all simple triods in $X$  is  a Borel subset of $\KK(X)$ and the mapping ${\bf c}:  \TT \to X$,   which assigns to any simple triod  $T$ its center,  is   Borel. 
 }\end{ssect}

   \begin{ssect}   \label{component}  
  Arc-components: {\rm    
   We denote by $E_X$  the arc-connection equivalence relation  on   $X$.  If $X$ is a Polish space then
$E_X$ is analytic as the projection on $X^2$ of the Borel set  

\ct{$B=\{ \bigl((x,y), J\bigr)\in X^2 \times \cp\JJ(X): \ e(J)=\{x,y\}\}.$}    
 
\nd  Any  arc-component $C$ in a space $X$ is:

-- either  a singleton,

-- or  admits  a one-to-one continuous parametrization   $\ph: I\to C$ where $I$ is a (closed, open, half-open)  interval in $\R$  or the unit circle,   and we shall then  say that $C$  is a {\it curve},

-- or else   contains a {\it triod}   and we shall then  say that $C$  is a {\it triodic  component}.

\sk

We  denote by $\Theta^X$ the union of all triodic  components of $X$, 
 to which we will refer    
 as  the {\it triodic part}; and the space $X$ is said to be   {\sl atriodic} if $\Th^X=\0$.
If   the  space $X$  is    Polish then the equivalence relation $E_X$ is analytic, hence  each triodic  component,  as well as  the triodic part,   is analytic. 
Note that    there exist  in $\R^3$ compact subsets with non Borel arc-components (see \cite{ks}). 
     
\begin{ssect} \label{can} 
Canonical arc-metrics and arc-topologies:
{\rm Given   any metric space  $(X,d)$ we can consider the  mapping $\d:X\times X\to [0,+\infty]$ defined by
$$\d(x,y)=\inf\{ \diam (H): \ H  \ \text{arc-connected s. t.}  \  \{x,y\} \subset H\subset X\}$$
where  $\inf \0=\infty$. So  if $x\not=y$  are in the same arc-component 
then $$d(x,y)\leq \d(x,y)=\inf\{ \diam (J): \ J \in   \JJ(X) \  \text{s.t.} \ e(J)= \{x,y\}\}<\infty$$
and   if  not  then $\d(x,y)=\infty$.
Moreover   setting by convention $\a+\infty=\infty$ for any $\a\in [0,\infty]$ we have for all $x,y,z\in X$
$$\d(x,z)\leq \d(x,y)+\d(y,z)$$

Hence strictly speaking $\d$  is not a distance on $X$, but $\d$ induces  a distance  on any arc-component   of $X$. We shall refer  to $\d$  as  the  {\it canonical arc-pseudo-metric defined by} $(X,d)$, and to its restriction  to any subset $A$ of some arc-component  of $X$   as  the  {\it canonical arc-metric  on $A$ defined by $d$}. 

The set of all $\d$-open balls with finite radius, constitutes  a basis  of  a metrizable topology $\tau$ on   $X$, finer than   the initial  topology $t$  defined by  $d$.  Note that: 
$$x=\lim_n x_n \ \text{in} \   (X,\tau)\iff \left\{
\begin{array}{l} \Ex (J_n)_n   \ \text{in} \  \JJ(X,t): \  
  \Vx n ,  \ e(J_n)=\{x,x_n\} \\ \text{and} \   \{x\}=\lim_n J_n \ \text{in} \ \KK(X,t)\end{array}\right. $$ 
so  the    topology $\tau$ can in fact  be defined directly from the   topology $t$, 
 and we shall say that $\tau$ is the  {\it canonical arc-topology} defined by $t$.

 We emphasize   that   even  if the initial topology $t$    is  separable  the topology $\tau$  is not  in general.
For example  if we fix in the unit   circle  in      $\R^2$ a  copy $C$ of the Cantor space  then the  union of all rays    joining the center to  an element   of $C$,  is an arc-connected  compact space  $X$.  But if $\d$  denotes the canonical  arc-metric on $X$ defined by the euclidean metric  then  for any two distinct  element  $a\not=b$ in $C$, $\d(a,b)\geq  1$.

The canonical pseudo-metric $\d$   was   introduced in \cite{dsr} for a Polish planar space $X$. But  as the reader can easily check the following    properties extracted from \cite{dsr}, that we state  without proof,   do not rely on this additional assumption.

\begin{thm}\label{delta} 
Let  $\d$   be the  canonical arc-pseudo-metric defined by the metric space $(X,d)$, and let $t$ and $\tau$ be respectively the $d$-topology  and $\d$-topology on $X$.
\begin{enumerate} 
\item  If  $(X,d)$ is complete  then $(X,\d)$ is complete
\item $\JJ(X,t)=\JJ(X,\tau)$,  
 \item For any arc     $J\subset X$, $d$-$\diam(J)$= $\d$-$\diam(J)$
\item For any arc   $J\subset X$,  $(J,t)= (J,\tau)$
\item All open $\d$-balls with finite radius are arc-connected,  
\end{enumerate}
\end{thm}
}\end{ssect}

Note that by property (2)    the  arc-connected subsets of  $(X,d)$ and $ (X,\d)$ are the same, so  the reference  to  arcs and arc-connectedness in  the following properties   is  non ambiguous. 
} 
  \end{ssect} 

\sk 

\nd{\it The   plane arc-connection relation:}  All specific properties of the arc-connection relation in the plane are due to the following fundamental property  of the plane topology.      

\begin{thm}\label{moore} {\sc (Moore)}  
Any family of pairwise disjoint triods in the plane is countable.
\end{thm}

 \nd In  particular   any planar  set admits at most countably many triodic components.

\begin{defin} \label{ker} 
  If $\tau$ is  the   canonical   arc-topology  of  some     space  $(X,t)$, 
the {\em triodic kernel} of $X$, that we denote by $\S^X$, is the $\tau$-closure of the set  of all centers of   triods in $(X,t)$.  
\end{defin}

The  following theorem is  a synthetic summary of the main results   of \cite{dsr}.

 \begin{thm} \label{dsr}  
 Let  $(X,t)$ be a planar Polish space, and let $\tau$ be the   corresponding arc-topology.

 a)  $\S^X$   is  $\tau$-separable, hence $(\S^X, \tau)$ is  a  Polish space, and 
$\S^X$  is a Borel subset of $(X,t)$.
  
b)  The set $B=\{(x,J)\in X\times \JJ(X): \ J\subset X, \ e(J)=\{x,y\}  \et  J\cap \S^X=\{y\}\}$
   is Borel and 
   
   \quad for all  $x\in X\mns \S^X$,  $\card(B(x))\leq 2$.   
   
  c) The equivalence relation $E_X$ is Borel.
 \end{thm}

\begin{rem} \label{psi}{\rm
Note that the projection on $X$ of the set $B$ in property b) is the set $\S'=\Th^X\mns \S^X$.  
Since $B$ is a Borel set with finite sections then  it admits a Borel uniformization. Hence there exist Borel mappings 
 $\Psi:\S'\to \JJ(X)$ and $\psi:\S'\to \S^X$ \st for all $x\in \S'$,  $e(\Psi(x))=\{x, \psi(x)\}$ and $\Psi(x)\cap \S^X=\{\psi(x)\}$.
}
\end{rem}

 \section{The partial operation   $\bigvee$ on  oriented arcs.}\label {sVee}

\begin{ssect}    
Oriented  arcs: {\rm  An {oriented arc} is a triple $\vec I=(I,a,b)$ where $I$ (the {\it domain} of $\vec I$) is an arc and  the two elements set $e(I)=\{a,b\}$  is ordered  by the pair  $(a,b)$.  
We shall  then say that $I$ is an {\em arc joining $a$ to $b$}   
and  set: 
$${dom}\, (\vec I)=I  \ ; \   e(\vec I)=(a,b)  \ ; \   a=e_0(\vec I) \et  b=e_1(\vec I)$$ 
so $e(I) $ is a two elements subset  of $X$, while $e(\vec I)\in X^2$. 
We will denote by $ \f$ the flip operation which assigns to any  oriented arc $\vec I = (I, a,b)$ the arc $\f(\vec I)=(I,b,a)$.

Given any oriented arc  $\vec I=(I,a,b)$    the     relation on $I$   defined by:
$$ x\leq _{\vec I} y\iff  \Ind I {a,x}\subset \Ind I {a,y}\iff  \Ind I {y,b}\subset \Ind I {x,b}$$
is a total ordering and we denote by $<_{\vec I}$ the corresponding strict ordering, to which we  refer as
{\em  the} total order  on $I$ defined by $\vec I$. 
In the sequel the notation $\vec I$ will always   suppose  implicitly that ${dom}\, (\vec I)=I$. 
Also for simplicity, when there is no ambiguity on the orientation  of  $I$
we shall set $e_0 (I) =e_0(\vec I)$, $e_1 (I) =e_1(\vec I)$, and  denote by  $<_I$  the  total  order on $I$ defined by $\vec I$. 

The set    $\vec\JJ(X)=\{(I,a,b)\in \JJ (X)\times X^2:  \ e(I)=\{a,b\} \}$ of all oriented arcs is  clearly Borel. Hence  
the set  $\{(\vec I, x, y)\in  \vec\JJ(X)\times X^2: \  x<_{\vec I}y  \}$ 
 is  Borel  too,  and the flip operation $ \f$ is Borel.
}
\end{ssect}

\begin{defin}  \label{defvee}
Given two   oriented arcs   $\vec  I_0 =(I_0,a_0,b_0)$ and $\vec   I_1 =(I_1, a_1,b_1)$   \st  $I_0\cap I_1\nv$  and $a_0\not=b_1$ we define the   oriented arc:
$$\vec   I_0  \vee \vec  I_1 =(I, a_0, b_1) \ \text{\rm with } \ 
 \left\{ 
\begin{array}{l}
I=I_0^{\{a_0, c\}}\cup I_1^{\{c, b_1\}}\\
\et  \\  c=\min_{<_{\vec I_0}} (I_0\cap I_1) 
\end{array}
\right. 
$$
 \end{defin} 

\begin{figure}[ht]
\begin{center}
\caption{}
\includegraphics[width=10cm]{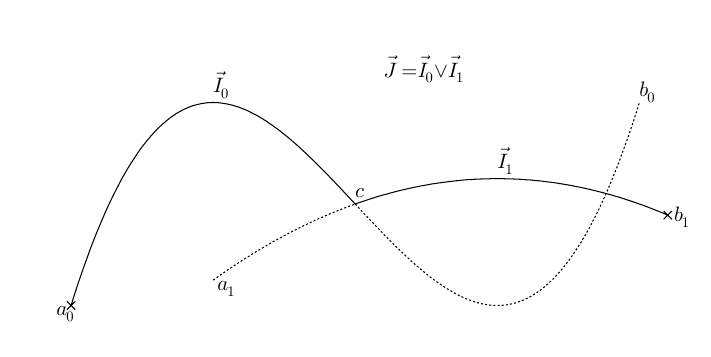}
\end{center}
\end{figure}
It follows readily from the previous definition   that if  $\vec   I_0  \vee \vec  I_1$ is   defined  then 
$$  \dom(\vec   I_0  \vee \vec  I_1) \subset    I_0  \cup  I_1  \ett  e( \vec   I_0  \vee \vec  I_1 )=(a_0,b_1)$$

\begin{defin} \label{defVee}  
a) For  any finite  sequence  $(\vec I_m)_{m\leq n}$  of   oriented arcs, we   define   inductively  the   oriented arcs:
    $$ \vec J_0=\bigvee_0\vec I_0  = \vec   I_0 \ett   \vec J_n = \bigvee_{m\leq n}  \vec   I_m=\vec J_{n-1} \vee \vec I_n.$$ 
 
 b)  If  $(\vec I_n)_{n\in \wo}$  is an infinite sequence in $\vec \JJ(X)$,  we shall say that $\vec J= \bigvee_{n\in \wo} \vec I_n $ is defined if   for all $n$,    $\vec J_n=\bigvee_{m \leq n} \vec I_m$ is  defined, and  $\vec J=\lim \vec J_n$.
\end{defin} 

\begin{rems} \label{remV}
{\rm a) We emphasize that  $\bigvee$ is a {\em partial operation}, so  $\vec J_n=\bigvee_{m\leq n}  \vec   I_m$ is defined only if all the terms in its  definition are defined, that is if for all $0< m\leq n:$
  $$ J_{m-1}\cap I_{m}\nv   \ett   e_1(\vec J_m)\not= e_0(\vec J_0)$$
and if so then:
$$ J_n\subset \bigcup _{m\leq n}   I_m \ \    \ett  e(\vec J_n)= (e_0(\vec I_0 ), e_1(\vec I_n)) $$

 \sk
 
b) If  for all $0< m\leq n$,    $ e_0(\vec I_{m})=e_1(\vec I_{m-1})\not = e_0(\vec J_0)$ then
$e_1(\vec I_{m-1})\not  = e_0(\vec J_m) = e_0(\vec J_0)$  and 
$ e_1(\vec J_{m-1})=e_0(\vec I_m)\in J_{m-1}\cap I_{m}$, hence $\vec J_n=\bigvee_{m\leq n}\vec I_m$ is defined.
}
\end{rems}

 \begin{lem} \label{BorV} 
 The set of all  finite or infinite sequences 
 $(\vec I_n)_{n<N}$ in $\vec \JJ(X)$ such that 
 $\vec J_N=\bigvee_{n<N} \vec I_n$ is   defined, is Borel,   and the mapping   which assigns 
 $\vec J_N$ to $(\vec I_n)_{n<N}$  is Borel.
 \end{lem}

 \begin{proof}
  It follows from its  definition   that the set $\DD_2$ of all pairs $(\vec I_0, \vec I_1)\in \bigl(\vec \JJ(X)\bigr)^2$ \st  $\vec   I_0  \vee \vec  I_1$ is   defined, is Borel.  Since the $\cap$ an $\cup$ operations  in $\KK(X)$, and the mapping $(I, a,b)\mapsto I^{\{a,b\}}$  are  Borel, then    the mapping $(\vec I_0, \vec I_1)\mapsto \vec I_0 \vee \vec I_1$ is Borel on $\DD_2$. 
This proves the lemma  for $N=2$, and the general finite case follows by a straightforward induction; and the 
infinite case follows from the definition.
 \end{proof}

 \begin{ssect} 
 Subdivisions:
 {\rm Let $\vec J=(J,a,b)$ be an oriented arc  and let $<$  be  the total order on  $J$ defined by $\vec J$.
A  {\em subdivison} of $\vec J$ is  a   finite sequence  $(\vec J^k)_{0\leq k\leq \ell}$ in $\vec \JJ(J)$ \st:   
$$a=e_0(J^0)<  e_1(J^0) = e_0(J^1) <  \cdots <   e_1(J^{k-1})=e_0(J^k)<    \cdots < e_1(J^{\ell-1})=
e_0(J^\ell)<  e_1(J^\ell)=b$$   
 }
 \end{ssect}

 \begin{thm}\label{Vee} 
 Let  $X$ be a  space and    $( \vec I_n)_{n\in\wo}$   be a  sequence   in $\vec \JJ(X)$ satisfying :
 
 a)  for all $n$,    $\vec J_n=\bigvee_{m \leq n} \vec I_m$ is  defined,

 b)  the sequence  $(I_n)_{n\in \wo}$ converges      in $\KK(X)$  to a singleton  $\{b\}$ with
 $b \not =a= e_0(I_0)$.

  \nd Then  $\vec J= \bigvee_{n\in \wo} \vec I_n $ is defined   and  $ e(\vec J)=(a,b)$.
\end{thm}

 \begin{proof}  
 Set for all $n$,    $\vec I_n=(I_n, a_n, b_n)$ and  $\vec J_n=(J_n,a, b_n)$, so   $b_n\not=a$.  
  For all $s, t\in \seq$,  we 
 denote by $s\wedge t$ the largest initial segment of $s\cap t$,  and if $s\nv$   we set:
$s^*=s_{\vert \abs s -1}$.

\mk
 
Then starting from $s_0=\<0\>$ and $\hat s_0=\<J_0\>=\<I_0\>$,  we construct inductively two 
 sequences $(s_n)_{n\in \wo}$ and $(\hat s_n)_{n\in \wo}$  in $\seq\mns\{\0\}$ and $\JJ(X)$ respectively   \st  for all $n$,  $\abs{\hat s_n}=\abs{s_n}$   and  $\hat s_n=(J^j_n)_{j<\abs{s_n}}$ is a subdivision of $\vec J_n$, 
and     for all $m<n$: 
\begin{enumerate}
\item   $s_n$ is an increasing sequence of integers $\leq n$ with     $s_n(0)=0$, 
 \item  for all $j <\abs{s_n}$, $ J^j_n\subset  I_{s_n(j)}$,  
\item $s_{n}\preceq s_{n-1}$ or  $s_{n}= s\cat \<n\>  $ with $s\preceq  s_{n-1}$,
\item  if    $\abs{s_m\wedge s_n}=k+1$  then  for all $j<k$, $J^j_{n}  = J^j_{m} $,   and   $ J^k_{n}  \subset  J^k_{m} $ 
         
\end{enumerate}
    \sk

    Suppose that   $(J^k_n)_{k<\abs{s_n}}$   is already defined satisfying conditions (1) to  (4). 
Then by definition $J_{n+1}=J_n^{\{a,c\}}\cup I_{n+1}^{\{c, b_{n+1}\}}$ with $c\in J_{n}\cap I_{n+1}$. Let
 $$k=\min\{-1\leq j<\abs{s_n}: \  c\,\leq_{J_n} e_1(J_{n}^{j})\}<\abs {s_n} \}$$ 
with the convention $  e_1(J_{n}^{-1})= e_0(J_{n}^0)=a$, so $k=-1$ if $c=a$. 
We then distinguish two cases:
 
 \sk
 
\ $(i)$  if $c=b_{n+1}$   then  
  $J_{n+1}=J_n^{\{a,b_{n+1}\}}=J_n^{\{a,a_{n}^{k}\}}\cup J_{n}^{\{a_{n}^{k},b_{n+1}\}} =
\bigl(\bigcup_{j<k} J_n^{j}\bigr)\cup J_{n}^{\{a_{n}^{k},b_{n+1}\}}  $  
 
 \sk
 
  and we  set:    
     $s_{n+1}= {s_{n}}_{\vert k +1}      \et   J^j_{n+1}  =\left\{ 
  \begin{array}{ll} 
 J^j_n & \si  j<k\\
 J_{n}^{\{a_n^k,b_{n+1}\}}  &  \si  j= k   \end{array} \right.$

\sk
 
$(ii)$  if  $c\not=b_{n+1}$ then
 $J_{n+1} =J_n^{\{a,c\}}\cup I_{n+1}^{\{c,b_{n+1}\}} =
\bigl(\bigcup_{j<k} J_n^{j}\bigr) \cup J_n^{\{a_{n}^{k},c\}} \cup I_{n+1}^{\{c,b_{n+1}\}} $ 

\sk 

and we  set:
   $ s_{n+1}= {s_{n}}_{\vert k +1}  \cat \<n+1\>   \et    J^j_{n+1} =\left\{ 
  \begin{array}{ll} 
   J^j_n& \si  j<k\\
  J_n^{\{a_n^k, c   \}} &  \si  j= k\\
    I_{n+1}^{\{c,b_{n+1}\}} &  \si  j= k + 1  \end{array} \right.
  $    

\sk

 In particular if $c=a$ then $k=-1$, so ${s_{n}}_{\vert k +1}=\0$ and  ${s_{n+1}} =\<n+1\>$ is of length 1,  and $J^0_{n+1} =I_{n+1}^{\{a,b_{n+1}\}}$ is the unique element of $\hat s_{n+1} $. This finishes the definition of 
 $s_{n+1}$ and  $\hat s_{n+1}$ which clearly satisfy conditions (1), (2), (3); and we now prove  condition (4).
  
So suppose that  $m<n+1$ and  $u=s_m\wedge s_{n+1}$, then by condition (3)  $u\preceq s_{n}$
hence $u\preceq s_m\wedge s_n$. If  $m<n$ then   (4)     follows from  the induction hypothesis; and   if $m=n$ then  $u\preceq s_n$ and    (4)   follows from the definition of $s_{n+1}$.

\mk

This ends up   the construction of the sequences $(s_n)_{n\in \wo}$ and $(\hat s_n)_{n\in \wo}$, and we set  $s_{-1}=\hat s_{-1}=\0$.

\begin{lem}  For all $n$:
\begin{enumerate}[resume]
\item if $\ell<m <n$ and   $s_\ell \preceq s_n$ then    $s_\ell \preceq s_m$                                                                
\item if $\ell< n$ and   $s_\ell \prec t\prec  s_n$ then there exists $m$  \st   $\ell<m <n$  and  $t= s_m$                                                                
\end{enumerate}         
\end{lem}

\begin{proof} The proof is by induction on $n$. For $n=-1$ the statements are trivial.
So suppose they are true for $n$.

\sk

 Suppose that  $\ell<m<n+1$  and $s_\ell \preceq s_{n+1}$. Since $\ell<n+1$ then necessarily 
 $s_\ell \preceq s^*_{n+1} \preceq s_{n}$.
 If $m<n$ then  by the induction hypothesis   $s_\ell \preceq s_{m}$; and 
if $m=n$  then $s_\ell \preceq s_n=s_m$. This proves (5) for $n+1$.

\sk
 
 Suppose that  $\ell< n+1$ and   $s_\ell \prec t\prec  s_{n+1}$. Since $t\not=s_{n+1}$ then 
 $t \preceq s^*_{n+1} \preceq s_{n}$. Then either $t =s_{n}$ and we are done; or else
 $t \prec s_{n}$ and then  by the induction hypothesis  there exists  $m$  \st   $\ell<m <n$  and  $t= s_m$. This proves (6) for $n+1$.
\end{proof}

\begin{lem}\label{III} 
One the following two alternatives holds:

\ (I) There exists $m\geq -1$ \st  the set   $N_m=\{n>m: \   s_n= s_m  \ou s_n^*= s_m \}$  is infinite,

 (II) There exists an increasing sequence $(n_i)_{i\geq _-1}$ such that for all $i$, and all $n\geq n_i$, 
$s_{n_i}\preceq s_{n}.$
\end{lem}

\begin{proof}  
Note that by  condition (6) the set $S=\{s_n;  \ n\geq -1\}$ is a tree, and consider the set 
$S'=\{s_m\in S:  \ \Vx n>m, \ s_n\succeq s_m\}$ which is nonempty since $\0=s_{-1}\in S'$.
It follows from condition (5) that if $\ell<m $ and   $s_\ell, s_m$ are in $S'$   then $s_\ell\preceq s_m$, hence
$S'$ is a $ \preceq$-chain.  

\sk

-- If $S'$ is finite, Let  $s$  be  the  $ \preceq$-maximum of $S'$ and 
$m=\min \{\ell: \  s_\ell=s\}$.      Then for all $n>m$, 
$s_m\preceq s_n$,  hence either $s_n=s_m$, or  $s_n^*=s_m$, or $\abs{s_n}> \abs{s_m} +1$ and then by condition (6) (applied to $t={s_n}_{\vert \abs{s_m} +1}$)   there exists   $n'$ such that $m<n'< n$ and $t=s_{n'}$ hence $s_{n'}^*=s_m$; and it follows that alternative (I) holds.

\sk

-- If $S'$ is infinite then there exists a unique increasing sequence $(n_i)_{i\in \wo}$ in $\wo$ such that
$S'\mns\0=\{s_{n_i}; \ i\in \wo\}$ and alternative (II) clearly holds.
\end{proof}

 For the rest of  the proof of Theorem \ref{Vee} we distinguish two cases according to Lemma \ref{III}. 
 Also for more clarity we shall split   alternative  (I) into two sub-alternatives.

\mk

{\it Case (I.a):} {\it There exists $m$ \st the set  $N=\{n \in \wo : \ s_n=s_m\}$ is infinite}

\sk
Set   $\abs{s_m}=k+1$; then  $H=\bigcup_{j<k} J^j_m$  is an arc with  
$e(H)=\{a, b^{k-1}_m \}$  and for all $n\in N$,  $J_n=H\cup J^k_n$, 
with  $e(J^k_n)= \{a_n^{k}, b^k_n\}= \{b^{k-1}_m, b^k_n\}$
 and  $ J^k_n\supset J^k_{n+1}$,  hence $J'=\bigcap_{n\geq m} J^k_n$ is an arc with $e(J')= \{b^{k-1}_m, b\}$ 
and  $J=\lim_{n\in N} J_n =H\cup J'$ exists and is  an arc with $e(J)= \{a, b\}$. 

 Let $(n_i)_{i \in \wo}$ be  the    increasing enumeration of $N$. 
 If $n_i\leq n <n_{i+1}$ then by  definition of $N$  we have   $s_{n_i} \preceq s_n $;  hence by condition (4) $J_{n_i} \subset J_{n}$ and    by condition (2) 
  $$J_{n} \mns J_{n_i} \subset \bigcup \{ I_{s_n(j)}: \ {\abs {s_{n_i}} \leq j< \abs {s_{n}}}\}\, .$$
  
  If  $k_i=\abs {s_{n_i}}-1$ then  $p_i=s_n(k_i)=s_{n_i}({k_i})\nearrow \infty$, and 
 since $\lim I_n$ is a singleton then
$ \diam (J_{n} \mns J_{n_i})\leq \diam (\bigcup_{p\geq p_i} I_{p})\searrow 0$.
Moreover  since $e_0( I_{p_i})=e_0( I_{s_{n_i(k_i)}})= e_0( J_{n_i}^{k_i})$ then  $I_{p_i}\cap  J_{n_i}\nv$.
Hence    if $d_H$ is the Hausdorff distance on $\KK(X)$ associated to some compatible distance on $X$, then  $d_H(J_n,J_{n_i}))= \diam (J_{n} \mns J_{n_i})\searrow 0$,
 which proves that $\lim_n J_{n} =\lim_i J_{n_i} =J$.

 \sk

{\it Case (I.b):} {\it There exists $m$ \st the set  $N^*=\{n \in \wo : \ s_n=s_m\cat \<n\>\}$ is infinite}

\sk

The argument is essentially the same. Setting $\abs{s_m}=k+1$ and  $H=\bigcup_{j<k} J^j_m$  as in     (I.a)
observe that  if $n\in N^*$, since   $s_n=s_m\cat \<n\>$,  then    $\abs{s_m}=k+2$ and  
   $J_n=H\cup J^k_n \cup J^{k+1}_n$. By the same arguments as in (I.a),  $J=\lim_{n\in N^*} H\cup J^k_n$ is an arc with $e(J)=\{a,b\}$. But   by condition (2) the additional term  $ J^{k+1}_n$ is a subset of $I_{s_n(k+1)}$ which by hypothesis b) of the theorem converges to the singleton $\{b\}$; and it follows that  $J=\lim_{n\in N^*} J_n$.
The rest of the argument is exactly the same as in (I.a)

\sk

{\it Case II:}  {\it  $S'$ is infinite.}

Set  $S'\mns\{\0\}=\{s_n;   n  \in M   \}$ and  for all   $n\in M$,  
   $\abs{s_{n}}=k_n+1$ and
$H_n=\bigcup_{j< k_n} J_{n}^j\subset J_{n}$. By  condition (4)  $(H_n)_{n\in M}$ is an increasing sequence 
of arcs with  $e(H_n)=\{a, c_n\}$  and  $c_n=b^{k_n-1}_{n}=  a^{k_n}_{n}\in I_{s_{n}(k_n)}$. Hence
  $\lim_{n\in M} c_n=b$ and $J=\adh{\bigcup_{n\in M} H_n}=\lim_{n\in M} H_n$  exists and is an arc with 
$e(J)=\{a,b\}$.  Finally    as   in Case (I.a)  if $(m_i)_{i\in \wo}$ is the increasing enumeration of $M$ then it follows from   condition (2) that for all $m=m_i\leq n<m_{i+1}$,  if $p_m=s_{m}(\abs{s_m}-1)$ then 
$d_H(J_{n},H_{m_i})=\diam(J_{n}\mns H_{m_i})\searrow 0$  which proves    that  $\lim_n J_{n}=J$.
 \end{proof}
 
\begin{rem} \label{notaV}{\rm 
Let $(I_n)_{n\in \wo}$ be  an   infinite sequence in  $\JJ(X)$   \st :

a)  for all $n$, 
   $e(I_n)=\{a_n,a_{n+1}\}$ for some sequence $(a_n)_{n\in \wo}$  in $X$. 
 
 b) the sequence  $(I_n)_{n\in \wo}$ converges      in $\KK(X)$  to a singleton  $\{b\}$ with
 $b \not =a_0$.
 
 \nd Then by Remark \ref{remV}.b)  for all $n$, $\bigvee_{m \leq n} \vec I_m$ is  defined, hence by
 Theorem \ref{Vee} the oriented arc  $ \bigvee_{n\in \wo}(I_n,a_n,a_{n+1}) =(J,a_0,b)$  is defined; and we shall write
$\bigvee_{n\in \wo}   I_n=J$.

Moreover by Proposition \ref{BorV} the mapping which to any such sequence $(I_n)_{n\in \wo}$ assigns the arc $\bigvee_{n\in \wo}   I_n$ is  Borel.
}
\end{rem}

 \section{ Arc-lifting}\label{Lift} 
 
 {We recall that $E_X$ denotes the arc-connection equivalence relation  on a given space $X$.}

 \begin{defin}  
 Given any   subset $S$    of some   space $X$, an  \emph{arc-lifting of $S$  in  $X$}
 is a   mapping  $\psi: S^2\cap E_X\to \cp\JJ(X)$ \st  for all $(x,y)\in  S^2\cap E_X$,  
 $e(\psi(x,y))=\{x,y\}$.
  If $S=X$ we shall say that $\psi$  is an  \emph{arc-lifting of    $X$.}
\end{defin}

    If $X$ is   Polish   then   the   equivalence relation   $E_X$ is analytic,  and  it  follows  then  from 
    Theorem~ \ref{YvN} that   $X$  admits a  $\SSS$-measurable  arc-lifting. But  as we shall see next the existence of a Borel arc-lifting for a space is a very strong assumption.

\begin{prop} 
If a Polish space $X$  admits a  Borel arc-lifting   then the equivalence relation $E_X$ is Borel.
\end{prop}

\begin{proof} 
Suppose that $\psi: E_X\to\cp \JJ(X)$ is a Borel arc-lifting. Since $E_X$ and $\cp\JJ(X)$ are   Borel then by classical results $\psi$ admits a Borel extension  $\td \psi: \DD \to \cp\JJ(X)$  to a Borel domain 
$E_X \subset\DD\subset X^2$ \st for all    $(x,y)\in \DD$ ,  $e(\td\psi(x,y))=\{x,y\}$; hence necessarily $E_X=\DD$, and so 
$E_X$  is Borel.
\end{proof}

We recall (see Theorem \ref{dsr})  that  if $X\subset \R^2$ then $E_X$ is Borel, but  there are compact spaces   in  $\R^3$    with non Borel arc-components, hence which do not admit  a  Borel arc-lifting.

\begin{defin} 
Given any   subset $S$  of some   space $X$,  and   any element     $a\in X$,  a  \emph{$\star$-arc-lifting   of $S$  in $X$, with summit $a$},   is  a  mapping  
$\ph: S  \to  \hat\JJ(X)$   \st  for   all  
$x\in S$, $e(\ph(x))=\{a,x\}$.   
\end{defin}

Note that   
 if $ S$   admits a   $\star$-arc-lifting   with summit $a$ then $a$ is not necessarily an element of  $S$, but $S\cup \{a\}$ is a subset of some  arc-component of $X$. 

\begin{rem}\label{star} {\rm Let $\Theta$ denote the  partial mapping  which to any pair $(J_0,J_1)\in \JJ(X)$ such that $e(J_0)=\{a,x\}$ and $e(J_1)=\{a,y\}$  with $x\not=y$,   associates the arc $J$ defined by $(J, x,y)=(J_0,x,a) \vee (J_1,a,y)$.  Then for any   $\star$-arc-lifting  $\ph:  S  \to \cp \JJ(X)$  of $S$ in $X$,  with summit $a$, and $x\not=y$ in $S$, the arc 
$$J= \psi(x,y)=\left\{ 
\begin{array}{ll}
\Theta(\ph(x), \ph(y)) & \ \text{if} \ x\not=a   \et y\not=a\\
\ph(y) & \ \text{if} \ x=a \\
\ph(x) &  \ \text{if} \ y=a 
\end{array}
\right.$$
is well defined and $e(J)=\{x,y\}$. The mapping  $\psi:S^2\to \JJ(X)$ thus defined is clearly 
   an arc-lifting  of $S$ in $X$.  Moreover  it follows from
 Lemma \ref{BorV} that  if $\ph$ is Borel then $\psi$  is Borel too.
We shall refer to $\psi$ as {\em the canonical arc-lifting} defined by $\ph$.} \end{rem}


\begin{thm}\label{Slift} 
Let  $(M, \d)$  be   a  complete,   arc-connected   and locally arc-connected,   metric space. Then  any separable subset $S$ of $M$  admits a  Borel arc-lifting in some separable subspace $N\supset S$.
 \end{thm}

\begin{proof} Fixing an element $a\in S$ we shall construct a  $\star$-arc-lifting  of summit $a$, of $S$.
Since  the set of  isolated points in $S$ is countable we may suppose that $S$ has no  isolated points. 
Also replacing $S$ by  $\adh S$ we may  suppose that $S$ is closed.

Since $S$ is separable we can fix a countable family   $\UU$ of open arc-connected subsets of $M$ whose restriction to $S$ form a basis for the topology of $S$. 
Then starting from $U_\0=X$  we can construct  a family $(U_s)_{s\in \seq}$  in  $\UU$,  
satisfying for all $s\in \seq$ and all $n\in \wo$:

\sk

\begin{enumerate}
\item  $ S\cap U_s\nv$ , $\diam (U_s)<  2^{-\abs s}    \et  \adh{U_{s\cat \<n\>}} \subset U_s$,
\item  $ S\cap U_s \subset   \bigcup_{n\in \wo}   U_{s\cat \<n\>}$.
\end{enumerate}

\sk 

Then   for all $\s\in \seq$,  $\bigcap_{n\in \wo} U_{\s_\vert n}$ is a singleton $\{\ph(\s)\}\subset S$ and
 the mapping $\psi:\wo^\wo \to  S$  thus defined is   onto, continuous and open, hence admits a Borel section $\ph: S\to  \wo^\wo$. 

  Set $a_\0=a$. Since $S$ has no isolated points,    then for all $s\in T\setminus \{\0\}$ we can fix an element $a_s\in S\cap U_s\mns\{a_{s_{\vert k}}; \ k<\abs  s\}$. 
So   if $s^*=s_{\vert \abs s -1}$ then  by construction $a_s\not= a_{s^*}$ and  $\{a_{s^*}, a_s\}\subset  U_s$; and since $U_s$ is arc-connected  we can  fix an  arc $I_{s}\subset U_s$ \st     $e(I_{s})=\{a_{s^*}, a_s\}$, hence  
 $\diam(I_s)<2^{-\abs s}$.
Then  for any   $\s\in \wo^\wo$, $a_\s= \lim_k a_{\s_{\vert k}}$ exists,  and if 
$a_\s\not=a$ then (see Remark \ref{notaV})  the  arc $\bigvee_{k\in \wo} I_{\s_{\vert k}}= (J_\s, a, a_\s)$ is defined, and    the mapping $\s\mapsto J_\s$ is Borel.  
 
 Hence  the mapping $\Phi:S \to \hat \JJ(X)$  defined by  $\Phi(x) =J_{\ph(x)}$ if $x\not=a$ and  $\Phi(a)=\{a\}$  is  a Borel $\star$-arc-lifting with summit $a$, of   $S$   in the separable space  $N=S\cup \bigcup_{s\in \seq} I_s$,  so  by    Remark \ref{star},  $S$  admits a  Borel arc-lifting in $N$. 
 in $X$. 
\end{proof} 

\begin{corol} \label{loc} 
Any  locally  connected   Polish  space  admits a  Borel arc-lifting.  
\end{corol}
\begin{proof} 
It is well known that any  locally   connected Polish space is locally arc-connected.  
Hence  $X$ admits   countably  many   arc-components $\{C_n; \ n\in \wo\}$ and 
each  $C_n$   is a clopen subset of $X$, so  a Polish space.   Hence  by  Theorem \ref{Slift}  each $C_n$ admits 
 a Borel arc-lifting $\Psi_n$    in $X$,  and  then   $\Psi=\bigcup_{n \in \omega}\Psi_n$  is clearly a Borel arc-lifting   of  $X$. 
\end{proof}

\begin{corol}\label{kerarc}
Let $(X,t)$  be a Polish space, and let $\tau$ be the arc-topology   defined by  $(X,t)$. Then any $\tau$-separable subset  $S$ of    some arc-component $C$  of $X$    admits a  Borel arc-lifting in $X$.  
\end{corol}
\begin{proof} 
Fix a complete  distance  $d$ on $X$ compatible with $t$. Then the topology $\tau$  can be defined by  the pseudo-arc-metric $\d$  defined by  $(X,d)$, and $\d$ is a metric on 
 $C$ of $X$. Then  $(C ,\d)$ is complete, arc-connected   and locally arc-connected,   metric space,   hence  by Theorem \ref{Slift} there exists a $\tau$ separable space $N\supset S$ and  a Borel  arc-lifting mapping   from $(C,\tau)$ to  $\hat\JJ(N,\tau)$. Then since $C$ is $t$-Borel and  the identity mapping from $\hat\JJ(N,\tau)$ to $\hat\JJ(N, t)$ is continuous, then $\Phi$ induces a   Borel arc-lifting of $S=(S,t)$  in $X$.
\end{proof}

 \begin{thm} \label{liftplan} 
Any planar Polish space $X$  admits a Borel arc-lifting.
\end{thm}
  
\begin{proof} Let   $d$ be a compatible complete distance on $X$. We denote  by $Y$ and $Z$  respectively  the  atriodic  part  and  the triodic part of  $X$. Since $X$ is a plane subset then $Z$ admits at most countably many arc-components $(Z_n)_{n\in N}$, and since 
 the equivalence relation $E_X$ is Borel  (Theorem \ref{dsr}.c)) then   all the $Z_n$'s as well as $Y$ are Borel subsets of $X$. Note that $$E_X=E_Y\cup \bigcup_{n\in N} E_{Z_n}=E_Y\cup \bigcup_{n\in N} Z_n^2$$ 
and we shall define    Borel arc-liftings  in $X$ separately   for  $Y$ and  for each $Z_n$ .

For $Y$ this is straightforward:  Since the set
$$B=\{(x,y, J)\in Y^2\times \hat\JJ(X): \ J\subset X \et e(J)=\{x,y\}\}$$
is Borel  and for all $(x,y)\in Y^2\cap E_X$ the section $B(x,y)=\{J_{(x,y)}\}$ is a singleton, then    the mapping
 $ (x,y)\to J_{(x,y)}$ is 
 a  Borel arc-lifting  of   $Y$ in $X$. 

\sk

We now fix $n$,  set $C=Z_n$ and $S=\S^X\cap Z_n$,
and  fix some element $a\in S$. To construct an arc-lifting of $C$ in $X$, it is enough to  
construct a  $\star$-arc-lifting $\Phi:C\to \JJ(X)$ of $C$ in $X$,  and again we shall define $\Phi$ separately   on   $S$ and $C\mns S$.  

\sk

By Theorem ~\ref{dsr}.a)  the set $S$ is $\tau$-separable, then by
Corollary~\ref{kerarc} we get an arc-lifting  $\Psi_S:S^2\to \JJ(X)$ of $S$ in $X$. In particular 
$\Phi_{S}= \Psi_S(a,\cdot): S \to \JJ(X)$ is a  $\star$-arc-lifting  with   summit $a$ of $S$ in $X$, 
and we  now  construct  a  $\star$-arc-lifting  $\Phi_{S'}$ with   summit $a$ of $S'=C\mns S$ in $X$.

\sk
By Remark \ref{psi}  we can fix two Borel mappings $\Psi: S'\to \JJ(X)$ and  $\psi: S'\to  S$ \st
 for all $x\in S'$, if $J= \Psi(x)$ then    $e(J)=\{x,\psi(x)\}$  and $J\cap S=\{\psi(x)\}$.
Then for all $x\in S'$,  $x\not =a$, and since $e(\Phi_S(\psi(x))=\{a,\psi(x)\}$  then
   $\psi(x)\in \Psi(x) \cap  \Phi_S(\psi(x))$,  hence  the oriented arc 
 $ (\Psi(x), x, \psi(x)) \vee (\Phi_S(\psi(x)), \psi(x), a)$  is defined, and is of the form  $( \Phi_{S'}(x), x,a)$. 
Then  by Lemma~\ref{BorV} the mapping  $\Phi_{S'}:S' \to \JJ(X)$ thus defined is  Borel and by construction $\Phi_{S'}$  is  a
  $\star$-arc-lifting  with   summit $a$ of $S'$ in $X$.
 \end{proof}

\section{Uniform arc-liftings} \label{Uniform}

 \begin{defin}   
 Given a space $X$ and $\GG\subset \FF(X)$  let  
  $E_\GG =\{(F,x,y)\in \GG\times X^2: \  (x,y)\in E_F\}.$ 
 A {\em uniform arc-lifting} of  $\GG$ is  a mapping 
 $\Psi:E_\GG \to \hat \JJ(X)$  \st for any $F\in  \GG $,    the partial mapping $\Psi(F, \cdot): E_F\to  \hat \JJ(X)$ is an arc-lifting of  $F$ (in $F$).
  \end{defin} 

As we shall see the complexity of the set $\Carc(X)$ is intimately related to the complexity of a potential uniform 
arc-lifting of $\Carc(X)$. Note that if $X$  is a Polish  space and $\GG$ is an analytic subset of $\FF(X)$ then  the set $E_\GG$  is analytic, hence  the set   $\BB=\{\bigl( (F,x,y), J\bigr)\in E_\GG\times \JJ(X): \  J\subset F \et e(J)=\{x,y\}\}$    is analytic too. It follows then from Theorem \ref{YvN} that  $\GG$ admits a  $\SSS$-measurable  uniform arc-lifting.
 If   moreover $X\subset \R^2$  then by Theorem \ref{liftplan}  any $F\in \GG$, admits   a  Borel arc-lifting, and it is natural to ask whether in this context,  $\GG$ admits a Borel    uniform arc-lifting. As we shall see if the set $\GG$ is rich  enough  then even the existence of a bianalytic    uniform arc-lifting is a very strong requirement.

\begin{prop} \label{atr}  
For any  Polish space $X$ and all $n\geq 1$, the set 
$$\Carc^{[n]}  (X)=\bigl\{F\in \FF(X):  \ \Vx x,y\in F, \ \Ex ^{\leq n}\, J\in \JJ(X): 
\  J\subset F \et e(J)=\{x,y\} \bigr\}$$
  admits   a  bianalytic    uniform arc-lifting.
\end{prop}  

\begin{proof}    
 Note that the set      
 $\BB=\{(F, x, y, J)\in \FF(X)\times X^2\times \JJ(X):   \  J\subset F \et e(J)=\{x,y\}\}$ 
 is Borel hence  by   Proposition \ref{unicity} the set
  $$
  \mathcal C=\{(F, x, y)\in \FF(X)\times X^2: \ \Ex ^{\leq n}\, J\in \JJ(X):  \   (F, x, y, J)\in \BB\}
  $$
 is $\ca$ and  $\BB\cap (\mathcal C\times \JJ(X))$ is the graph of  a bianalytic mapping   
 $\Phi$ on $\mathcal C$, which to any $F, x, y)\in \C$ assigns a nonempty finite set
 $\{J_1, \cdots, J_n\}$ in $\JJ(X)$  with graph contained in $\BB$. Then  any   bianalytic  section $\Psi: \mathcal C \to  \JJ(X)$ of $\Phi$, given by   Proposition \ref{unicity},   is a  uniform  arc-lifting for $\Carc^{[n]}  (X)$.
 \end{proof}

\begin{prop} \label{GG} 
Suppose that  $\GG\subset \Carc(X)$ is in $\ca(\Carc(X))$, that is $\GG=\HH\cap \Carc(X)$ where
 $\HH$ is $\ca$. If $\GG$    admits a  bianalytic    uniform arc-lifting then   $\GG$ is $\ca$.
\end{prop}  

\begin{proof}  
Let    $\Psi: E_\GG\to \hat \JJ(X)$  be   a  bianalytic    uniform arc-lifting.
Since the condition ``$J\subset F \et e(J)=\{x,y\}$"  is Borel, then by Proposition \ref{ext} $\Psi$ admits an extension $\td \Psi$ to   a bianalytic   mapping with  $\ca$  domain $\DD\subset \HH \times X^2$ and     \st if  $J=\td \Psi(F,x,y)$ then $J\subset F$ and  $e(J)=\{x,y\}.$ 
Then the set 
  $\mathcal C_*=\{F\in \HH: \ \Vx x,y\in F, \ (F,x,y)\in \DD\}$ 
 is  $\ca$. If  $F\in \GG$ then $F\in \HH$  and $E_F=F^2$ hence  $\{F\}\times F^2\subset E_\GG\cap ( \HH \times X^2) \subset \DD$; so  
 $F\in \mathcal C_*$. Conversely if $F\in \mathcal C_*$ then then $F\in \HH$ and  the mapping $\td \Psi$ witnesses that  $F\in \Carc(X)$, hence $F\in \GG$. It follows that $\GG =\mathcal C_*$ is $\ca$.
 \end{proof}
 
   \begin{corol} \label{Cua}  
   For any  Polish space $X$, if $\Carc (X)$ admits  a  bianalytic    uniform arc-lifting  then $\Carc (X)$ is $\ca$.
\end{corol}  

\begin{corol} \label{CTh}  
For any  Polish space $X$, each of the following subsets of $\Carc(X)$
 is $\ca$ and admits a  bianalytic    uniform arc-lifting:

a) the set $\Carc^{\check \Th} (X)$ of all atriodic closed arc-connected subsets of $X$,

b) the set $\Carc^{0} (X)$ of all closed  arc-connected subsets of $X$   that  contain   no Jordan curve.

\nd In particular if $X$ is atriodic  or   contains no Jordan curve then $\Carc(X)$ is $\ca$.
\end{corol}

\begin{proof} 
a)  If $F\in \Carc^{\check \Th} (X)$  then any arc-component of $F$ is a curve (see Section \ref{component}) hence $\Carc^{\check \Th} (X)\subset \Carc^{[2]}  (X)$  so by  Proposition \ref{atr}, 
$\Carc^{\check \Th} (X)$ admits a  bianalytic    uniform arc-lifting. Moreover  $\Carc^{\check \Th} (X)=\HH\cap \Carc (X)$   where $\HH$ is the $\ca$ set of all atriodic closed subsets of $X$, hence by  Proposition  \ref{GG},  
$\Carc^{\check \Th} (X)$is $\ca$. 

b) The argument is  similar: if  $F\in \Carc^{0} (X)$  then any arc-component of $F$ is  the continuous image of subinterval of the real line, hence $\Carc^{0} (X)\subset \Carc^{[1]}  (X)$  and by  Proposition \ref{atr}, 
$\Carc^{\check \Th} (X)$ admits a  bianalytic    uniform arc-lifting. Also $\Carc^{0} (X)=\HH\cap \Carc (X)$   where 
$\HH$ is the $\ca$ set of all  closed subsets of $X$ not containing  any Jordan curve, hence by  Proposition  \ref{GG},  $\Carc^{\check \Th} (X)$ is $\ca$. 
\end{proof} 

\section
 {Uniform arc-liftings in the codes} \label{unifcode}  
 
 As we mentioned before   it  was  already known from the work of   Becker
 that the set  $\Carc(\I^2)$ is not $\ca$, hence  by Proposition \ref{GG},  $\Carc(\I^2)$ does not  admit  a  bianalytic    uniform arc-lifting. 
 Nevertheless we shall prove that for any Polish  $X\subset \R^2$  the set $\Carc(X)$ admits  a  weak  form  of bianalytic   uniform arc-lifting. 
Notice that since  by Corollary \ref{CTh}  the set $\Carc^{\check \Th} (X)$ does admit a   bianalytic    uniform arc-lifting    we can restrict the study  to the set  $\Carc^\Th(X)=  \Carc(X)\mns \Carc^{\check \Th} (X)$   of  all triodic arc-connected closed subsets of $X$. To state the main result precisely
we need to introduce  some preliminary notions.

    \begin{ssect} Coding  $\d$-separable subsets: 
 {\rm
 For all $F\in \FF(X)$ let  $d_F$ denote the distance on $F$ induced by $d$, and   let $\d^F$  be   the 
 arc-pseudo-metric defined on $F$ by the metric structure $(F,d_F)$. 
 Our goal  is to code pairs $(S, F)$ where $F\in \FF(X)$ and $S$ is a $\d^F$-separable closed subset  of $F$ contained in some arc-component     of $F$ (so that  $\d^F$ induces  a genuine distance on $S$). 
 Since the space  $(S, \d^F)$ is  separable,  it is entirely determined by  the  restriction of   $\d^F$  to any countable dense subset $D$  of $S$.
 }
  \end{ssect}

 {\it From now on we fix  a Polish space $X$ and set  $\FF=\FF(X)$, $\JJ=\JJ(X)$ and  $\A=   X^\wo  \times \JJ^\wo $.}
 
 \sk 

 \begin{defin}\label{code}  
  Let   $\C$ be the set of all $(\a,F)\in \A\times \FF$ with
  $\a=(\a_0, \a_1)=\bigl((a_n)_{n\in \wo},   (I_p)_{p\in \wo}\bigr)$   
  satisfying for all  $m,n\in \wo$: 
$$ \left\{
\begin{array} {l}
(1) \  a_n \in F \ \et \  I_n\subset F, \\
(2)  \  \d^F{(a_m,a_n)}=\inf\{\diam(I_p):  \ e(I_p)=\{a_m, a_n\}\}  \ \text{if} \     a_m\not=a_n.  
\end{array}
\right.
$$
We then   set:  $\a^*=a_0, \  D_\a=\{a_n; \ n \in \wo\}$,  $S_{(\a, F)}=\adh{D_\a}^{\d^F}$ and 
$S'_{(\a, F)}=F\mns \SaF$.
\end{defin}
\sk

One should  think to an element $(\a,F)\in \C$  as a code  for   the  separable metric  space $(S_{(\a, F)}, \d^F)$.  

\begin{lem} 
The set $\C$   is   $\ca$.
\end{lem}

\begin{proof} 
Condition (1)  is clearly  Borel. If  $\s(\a)$ denotes the  righthand side of the equation in   (2), then 
the  mapping $\a\mapsto \s(\a)$   is   Borel too;  and   it follows from the definition of $\d^F$  that  $\d^F_\a\leq \s(\a)$, hence     condition (2) is equivalent to:  
 $$ \Vx m, n \in \wo, \  \Vx J\in \JJ, \ \text{if} \ (J\subset F \et e(J)=\{a_m, a_n\}) \ \text{then} \
\s(\a) \leq \diam (J)$$
which is clearly $\ca$.
\end{proof}

{\nd Set:

\ct{$ \boldsymbol R= \{(\a, F, x)\in\C\times    X: \ x\in F\} \ , \
\SS= \{(\a, F, x)\in \C\times    X: \ x\in \SaF\} \et \SS'=\R \mns\SS$ } 

\nd so 

\ct{ $\R=\SS\cup \SS'$}
 
 \nd and we also set:
 
\ct{$ \boldsymbol R^{(2)}= \{(\a, F, x,y)\in\C\times    X^2: \ x,y \in F\} \et 
\SS^{(2)}= \{(\a, F, x,y)\in \C\times    X^2: \ x, y\in \SaF\}$}

\begin{lem}\label{RS2} For all $r>0$:

a)  The set  $\{(F, x,y)\in \FF\times X^{(2)}: \   x,y\in F \et \d^F(x,y) < r\}$ is  $\ana$.   

b) The set  $\{(\a, F, x,y)\in\SS^{(2)}: \   \d^F(x,y) < r\}$   is  bianalytic in $\SS^{(2)}$.
 
 \nd Similarly if we replace  $<$ by \;$\leq$.
 \end{lem}

\begin{proof} a) follows readily from the  definition of $\d^F$. For b) observe that  by the density of $D_\a$ in
$\SaF$ and the triangle inequality,  for $(\a, F, x,y)\in\SS^{(2)}$:

$\d^F(x,y)\geq r \iff \Vx r'<r, \  \d^F(x,y) > r'$

\hskip 20mm $\iff \Vx r'<r, \  \Ex \e>0,  m,n \in \wo, \
\left\{ \begin{array}{l}    \d^F(a_m,a_n) >r'+2\e,   \\
   \d^F(x, a_m) <\e \et  \d^F(y,a_n) <\e \, .
   \end{array} \right.$ 
\end{proof}

\begin{lem}  \label{RS}
a) The set  $\R$ is  a Borel subset of $\C\times   X$, hence a $\ca$ set.

b)  The set  $\SS$   is    in $\ana(\R)$, hence $\SS'$ is a $\ca$ set.
\end{lem}

\begin{proof}  Part a) is   obvious, and part b) follows from Lemma \ref{RS2} since  for $(\a,F,x)\in\R$:

 \ct{$(\a, F, x)\in \SS \iff    \Vx \e>0, \  \Ex n, \ \d^F(a_n, x) <  \e$.}  \end{proof}

\begin{thm} \label{Phi} 
There exists   a bianalytic mapping  $\Phi=(\Phi^{(0)},\Phi^{(1)}): \SS \to \wo^\wo\times \hat \JJ$  \st for all $(\a,F,x)\in \SS$ 
if $\Phi(x)= (\s,J)$ then:

1) $x=\d^F\text{-}\lim_n \a_0(\s(n))$

2) $J\subset F$ and $e(J)=\{\a^*, x\}$

\nd In particular the partial mapping $\Phi^{(1)}(\a,F, \cdot):  \SaF \to\hat \JJ$ is a Borel $\star$-arc-lifting  for $\SaF$ in $F$, with summit $\a^*$.
%
\end{thm}

  \begin{proof}  
 The argument is   a reminiscence of the proof of Theorem \ref{Slift}.
 We fix  a  Borel    bijection   $\rho: n\mapsto  (\ell_n, r_n)$ from $\wo$ onto $\wo\times\Q^+$. 
Let $\tau:n\mapsto \ell_n$ denote the first coordinate of $\rho$ and  for $(\a,F)\in \C$  set  for all $n$:  
 $a_n= \a_0(n)$ and  $c_n =\a_0(\ell_n) =\a_0(\tau (n))$.
 
 \sk
 
We then  define a tree $T=T(\a,F)\subset \seq$ on $\wo$  
as follows:

\sk

  $
\begin{array}{cl} 
 (i) & \quad  \ T\cap \wo^{1}= \{ \<m\>; \ m\in \wo\} \\ 
(ii) & \left\{\begin{array}{l}
\text{if} \  t=s\cat\<m\>\in T  \ \text{then:} \ \\
  {t\cat\<n\>\in T \iff  
    r_n< \min\{ r_m- \d^F  (c_n,c_m)  \ , \ \dfrac {r_m} 2\ , \  \{\d^F(c_n, c_{t(j)} );  \  j<\abs t\}} \end{array}\right.
\end{array}$
 
\sk
 For all $n$ let  $B_n=B_n(\a, F)$ and $\td B_n=\td B_n(\a, F)$ denote respectively the open $\d^F$-ball, and the closed $\d^F$-ball  in $\SaF$, of center $c_n$ and radius $r_n$.  Then it follows from   condition $(ii)$  
   and the triangle  inequality, that if $u=t\cat\<n\>=s\cat\<m\>\cat\<n\>\in T$    
  then $\td B_n\subset B_m$, $r_{n} < \dfrac {r_{m}} {2}$, and  
for  all $j<\abs t$, $c_n\not=c_{t(j)} $,  since $r_n< \d^F (c_n, c_{t(j)})$.

Hence   if  $s=\<n_0, n_1, \cdots, n_k\>\in T $  with $k\geq 1$ then 
 $B_{n_0}\supset  \tilde B_{n_1} \supset  B_{n_1}  \cdots\supset \tilde B_{n_k} \supset B_{n_k}$ and 
$r_{n_k} <  {r_{n_0}}\, {2^{-k}}$.  
Then  setting $S=S_{(\a,F)}$, for any $x\in S\cap B_{n_k}$ we can find $r\in \mathbb Q^+$ \st 
$$ r <\min\{ r_{n_k}- \d^F (x,\a_0{(n_k}))  \ , \  \dfrac{r_{n_k}} 2\ , \  \d^F{(x, \a_0({t(j)}))}  \ \text{for all} \  j<\abs t\}.$$
Also by  the density of $D_\a$ in $S$ we can find some $\ell$ \st the same inequality holds when replacing $x$ by 
  $a_p$;  and 
if  $(\ell,r)=(\ell_n, r_n)$  then  $s\cat \<n\>\in T$. Hence any $s\in T$ has a (strict) extension in $T$ and  
$B_{n_k}\subset  \bigcup\{B_n: \  s\cat\<n\>\in T\}$.     Since by condition $(i)$  
$S=\bigcup\{B_n: \ \<n\>\in T\}$  then   by induction
$S= \bigcup_{\s \in  \lceil  T\rceil }\bigcap_{j \in \wo}  B_{\s(j)} $ 
with $\bigcap_{j \in \wo} \downarrow B_{\s(j)}= \{b_\s\} $ for all $\s\in \lceil  T\rceil$.

So  for any  $x\in S$  there exists  at least  one  infinite branch $\s \in  \lceil  T\rceil $ \st $x\in  \bigcap_j B_{\s(j)}$, hence  $x=b_\s=\d^F$-$\lim_k c_{\s(k)}=\d^F$-$\lim_k a_{(\ell_{\s(k)})}$.   If  $\Phi^{(0)}(\a,F,x)=\s$ is  the lexicographical minimum of the set $\{\s'\in  \lceil  T\rceil : \ x=b_{\s'}\}$ then:
$$\begin{array}{ll}
\s(k)=n &\iff 
  x \in  \bigcap_{j \leq k}B_{\s^x(j)} \mns \bigcup\{B_{m}: \  m<n \ \text{and} \  \s_{\vert k}\cat \<m\>\in T\}\\
&\iff   \Vx j<k,  \  x\in B_{\s^x(j)}, \et ( \Vx m<n, \ \s_{\vert k}\cat \<m\>\not \in T \ou x\not\in B_m).
\end{array}$$
Since  by condition (2) of Definition \ref{code}, for all $m,n$,   the mapping  $F\mapsto \d^F(a_n, a_m)$ is Borel on $\C$, 
 then  for all $s\in \seq$,    the set
$\{ (\a, F)\in \C: \ s\in T (\a, F)\}$ is a  Borel  subset of  $\C$, hence the 
mapping $(\a, F)\mapsto T(\a, F)$ is  Borel on $\C$. Moreover  for all $n$, since $a_n\in \SaF$, then
by Lemma~\ref{RS2}  the set  $\{(\a, F, x)\in \SS: \  x\in B_n\}$ is  
   bianalytic  in $\SS$. Hence   the mapping $\Phi^{(0)}$ 
    is  bianalytic   on $\SS$.

 \sk  
 
 Then by condition (2) of Definition \ref{code} for all $(m,n)\in\wo^2$ we can define on $\C$  a   Borel  mapping $\mu_{(m,n)}:\C \to \wo$ by 
   $  \mu_{(m,n)}(\a,F)=\min\{ p: \      e(I_p)=\{a_m, a_n\} \et
\diam(I_p)< 2\,  \d^F(a_m, a_n,)\}  $.  If $(\a,F, x)\in \SS$  with
$\ph (\a,F, x)= \tau\rond \Phi^{(0)}(\a,F,x)= (\ell_{k_n)})_n \in \wo^\wo$  
 and  $p_n = \mu_{(\ell_{k_n} ,  \ell_{k_{n+1}})}(\a,F ) $ we set  $\psi  (\a,F, x)=  (p_n)_n $;     then the mapping $\psi: \SS\to  \wo^\wo$     thus  defined  is   bianalytic    on $\SS$. Moreover  by the definition of $\mu$, 
 for all $n$, $e(I_{p_n})=\{c_{k_n}, c_{k_{n+1}}\}$, \  
$\diam(I_{p_n})< 2\,  \d^F(c_{k_n}, c_{k_{n+1}})$ and    $x=\d^F$-$\lim_n c_{k_n}$, hence $x=d$-$\lim_n c_{k_n}$. 
If $x\not =a_0$   it follows  from  Remark \ref{remV}.\,b) that the sequence $(I_{p_n})_n$ satisfies the hypothesis of     
 Theorem \ref{Vee}, hence  the arc  $J=\Phi^{(1)}(\a,F,x)=\bigvee_{n\in \wo} I_{(p_n)}$ is well defined and satisfies $J\subset F$ and $e(J)=\{a_0,x\}$. Finally   by Lemma \ref{BorV},  setting $\Phi^{(1)} (\a,F, a_0)=\{a_0\}$, the mapping  
 $\Phi^{(1)}$  thus defined on $\SS$  satisfies the conclusion of Theorem \ref{main}. \end{proof}

\begin{corol} \label{SS'} The set  $\SS$ is $\ca$, hence  $\SS$ and $\SS'$ are  bianalytic subsets of  $\R$. 
\end {corol}

\begin{proof}  
Let    $\ph=\Phi^{(0)}: \SS\to \wo^\wo$ be the bianalytic mapping given by Theorem \ref{Phi}.
Since $\SS\subset \R$ and  $ \R$  is $\ca$, then $\ph$ admits an extension to a bianalytic mapping 
$\td\ph:   \Uu\to \wo^\wo$ with $\ca$ domain $\SS\subset   \Uu \subset  \R$.  
 Moreover observe that since $\d^F$ is complete on $F$ then for  $(\a,F,x)\in \R$ and   $\s\in \wo^\wo$: 
$$x=\d^F\text{-}\lim_n \a_0(\s(n)) \iff   x=\lim_n \a_0(\s(n)) \et  (\a_0(\s(n))_n \ \text{is $\d^F$-Cauchy}$$
and  it follows from Lemma \ref{RS2} b) that the  lefthand side of this equivalence is a bianalytic condition on $(\a, F, x)$.  Hence the set $\td\SS=\{(\a,F,x)\in   \Uu: \    x=\d^F\text{-}\lim_n \a_0(\td \ph((\a,F,x))(n))\}$ is $\ca$ and 
$\SS\subset \td \SS$; and since the converse inclusion is obviously true then 
 $\SS= \td \SS$,  and so $\SS$ is $\ca$. \end{proof}

From now on we restrict our study to closed subsets  $F$ of the given space $X$, which are: {\em triodic,  arc-connected, and $\d^F$-separable}.  For this we consider 
$$\C_\star=\{(\a, F)\in \C: \  F\in \Carc^\Th \ett  \SaF\supset \S^F\}$$
and for any set $\Q\subset \C  \times Y $ defined previously we set $\Q_\star=\Q\cap( \C_\star  \times Y )$.

\begin{thm} \label{Phi'} 
There exists   a bianalytic mapping  $\Phi':  \SS'_{\star} \to X\times \hat  \JJ$  \st for all 
$(\a,F,x)\in \SS'_{\star}$,   if $\Phi'(\a,F,x)=(y,J)$ then $y\in \SaF$, $J\subset F$ and  $e(J)=\{x,y\}$.
\end{thm}

\begin{proof} Set $\td \C_\star=\C_\star \times X^2 \times \JJ$ and
consider the sets:

\sk

 $\B=\{(\a, F, x, y, J)\in \td \C_\star :  \  J\subset F, \ e(J)=\{x,y\}, \ x\in \SaF'\et   y\in  J\cap \SaF\}$
   
$\B_0=\{(\a, F, x, y, J)\in \B :  \    \{y\}= J\cap \SaF\}$

$\B'_\e=\{(\a, F, x, y, J, z)\in \B  \times X: \   
z\in J\cap \SaF \et d(y,z)\geq \e\}$,
\ \  for  $\e>0$.
\sk

 By Corollary \ref{SS'} the sets $\B$ and $\B'_\e$ are bianalytic subsets of  
 $\td\C_\star$ and $\td \C_\star  \times  X$ respectively. 
 Moreover  by  Theorem \ref{delta} 4) for all 
 $(\a, F)\in \C$ and $J\in \JJ$ the set $J\cap \SaF$  is compact,
  hence  for all $(\a, F, x, y, J)\in \B$ the set $\B'_\e(\a, F, x, y, J)$ is compact too. It follows then 
from  Theorem \ref{unicity} that $\pi(\B'_\e)$ -- the projection of $\B'_\e$ on $\B$ -- is  bianalytic in $ \B$, hence in $\td\C_\star $. 
So $\B_0=\bigcap_\e \B \mns \pi(\B'_\e)$   is    bianalytic in $\B$.  
 
   \sk
 
Then  for all $(\a, F, x)\in \SS'_{\star}$,  since $F\in \Carc(X)$  and $\S^F\nv$,  we can join $x$   to $\S^F$, hence to $\SaF$, by an arc $J$,  and if $J$ is minimal for $\subset$ and $e(J)=\{x,y\}$  
then $(\a, F, x, y,J)\in \B$; hence  $\B(\a, F, x)$  is nonempty.
Observe also that if $(y, J)$ and $(y',J')$ are two distinct elements of $\B(\a, F, x)$  then  necessarily $J\not= J'$. 
Moreover  if $(\a, F, x)\in \SS'_{\star}$ then $x\in \SaF'$, and since $(\a,F)\in \C_\star$  then $x\not\in \S^F$, in particular $x$ is not the center of a triod.  It follows that for all $(\a, F, x)\in \SS'_{\star}$  the set $\B(\a, F, x)$  contains at least one element, and at most two  elements,  hence by   Theorem \ref{unicity} the set $\B$ admits a bianalytic  selection $\Phi'$ which satisfies the conclusion of Theorem \ref{Phi'}.
 \end{proof}

Recall that by definition 
$$(\a, F, x,y)\in \R^{(2)}_{\star} \iff  ( \a, F)\in \C_{\star}(X) \et (x,y)\in F^2$$
and

\ct{$( \a, F)\in \C_{\star}\iff \left\{
\begin{tabular}{l}
 -- $F$ is a   {\em triodic, arc-connected and   $\d^F$-separable} subset of $X$, \\
 --  $(\a,F)\in \C$ is a  code    satisfying  $\S^F\subset \SaF\subset  F$.
\end{tabular}\right.$} 

\begin{thm} \label{main} 
For any   Polish space $X$,  there exists a bianalytic mapping  $\Psi: \R^{(2)}_{\star} \to\hat  \JJ(X)$   \st   
for all  $( \a, F)\in \C_{\star}$,  the partial mapping $\Psi(\a,F,   \cdot): F^2\to \hat \JJ$ is a 
Borel arc-lifting for $F$.
\end{thm}
\begin{proof} Let $\Phi$ and $\Phi'$ be as in Theorem \ref{Phi} and Theorem \ref{Phi'}, and let 
 $\phi_0$  denote  the restriction of the mapping $\Phi^{(1)}$ to $\SS_\star$. So  for all $( \a, F)\in \C_\star$, the  mapping  $\phi_0  (\a,F, \cdot): \SaF \to \hat \JJ$   is  a $\star$-arc-lifting of $\SaF$ in $F$, with summit $\a^*$ .

Also for all   $( \a, F,x)\in \SS'_\star$  if  $y= \Phi'^{(0)} ( \a, F,x) $  and $J=\phi'( \a, F,x)$ is the arc   defined by
  $(J, x,\a^*)=\bigl (\Phi'^{(1)}( \a, F,x) , x,y\bigr)\vee (\Phi^{(1)}( \a, F,y), y, \a^*)$, 
then the mapping $\phi_1: \SS'_\star \to \hat \JJ$ thus defined is  
 bianalytic and for all $( \a, F)\in \C_\star$, the  mapping $\phi_1 (\a,F, \cdot): \SaF' \to \hat \JJ$ is a $\star$-arc-lifting of $\SaF'$ in $F$, with summit $\a^*$ .

  Since $\SS_\star$ and $\SS'_\star$  are   bianalytic subsets of $\R_\star$,   the  mapping  $\phi =\phi_0\cup \phi_1$ on  $\R_\star$ is bianalytic  too, and for all  $( \a, F)\in \C_\star$ 
 the  mapping $\phi  (\a,F, \cdot): F \to \hat \JJ$   is  a $\star$-arc-lifting of $F$, with summit $\a^*$; and since
 $F$ is Polish then $\phi  (\a,F, \cdot)$ is Borel.   Then if $\Psi(\a,F, \cdot, \cdot)$ denotes 
  the canonical  arc-lifting on $F$ defined by  $\phi (\a,F, \cdot)$ (see Remark \ref{star})   
  the mapping $\Psi: \R_\star \to\hat  \JJ(X)$   satisfies clearly the conclusion of Theorem \ref{main}. 
\end{proof}

\begin{thm} \label{sG} 
Let $X$ be   Polish space and suppose that  for any $F\in \Carc^\Th(X)$ the set $\S^F$ is $\d^F$-separable, and let
$\G\supset \ca$ be a given class.

a)  If  there exists a $\G$-measurable mapping ${\boldsymbol\s}: \Carc^\Th(X) \to \A$  \st for all $F\in \Carc^\Th(X)$,  $({\boldsymbol\s}(F),F) \in  \C_\star$  then
 the set $\Carc(X)$  is  in $\check{\mathcal A}(\G)$.   

b) If   there exists a $\G$-measurable mapping   ${\boldsymbol\s_0}: \Carc^\Th(X) \to X^\wo$ \st for all $F\in \Carc^\Th(X)$,  the $\d^F$-closure of the range of the sequence ${\boldsymbol\s_0}(F)$  contains the set  $\S^F$, then the set $\Carc(X)$  is  in  $\check {\mathcal A}(\check{\mathcal A}(\G))$.


\end{thm}

\begin{proof} Again  since by  Proposition \ref{atr}   the set $\Carc^{\check \Th} (X)$  is $\ca$,  to prove that
$\Carc(X)$ is in some class $\boldsymbol\Lambda \supset \ca$, we  only need to prove that the set  $\Carc^\Th(X)$  is  in $\boldsymbol\Lambda$.

\sk

a) Let $\Psi$ the mapping given by Theorem \ref{main}. Since the domain $\R_\star^{(2)}$ of the mapping $\Psi$ is a subset of the $\ca$ set $\R^{(2)}$,  
 and the  graph of $\Psi$ is a subset of the $\ca$ set  $$\Uu=\{(\a,F, x,y,J)\in \R^{(2)}: \   \S^F\subset \SaF, \  J\subset F \et e(J)=\{x,y\}\}$$
  we can extend $\Psi$ to  a bianalytic mapping  $\td \Psi: \Q\to \hat \JJ$ with   
  $\ca$ domain  $ \R_{\star}^{(2)}\subset \Q\subset  \R^{(2)}$ and 
  graph  also contained in $\Uu$. It follows that  the set
$$\D=\{(\a, F)\in \C: \   \Vx (x,y)\in F^2, \  (\a,F, x,y)\in \Q 
 \}$$
is $\ca$.    

For all $F\in \Carc^\Th(X)$, since  $({\boldsymbol\s}(F),F) \in \C_{\star}$,   then    
 by  Theorem \ref{main}   for all $x,y\in F$,  $(\a, F,x,y)\in \R_\star^{(2)}\subset \Q$; hence
$({\boldsymbol\s}(F),F) \in \D$.
Conversely if  $({\boldsymbol\s}(F),F) \in \D$ then by the choice of $\td \Psi$,  for all $x,y\in F$,  $\td \Psi(\a, F,x,y)$ is a sub-arc of $F$ with endpoints
$\{x,y\}$; hence $F\in \Carc^\Th(X)$. This proves that
$\Carc^\Th(X)=({\boldsymbol\s}, Id)\m (\D)$ and the conclusion follows from Proposition \ref{A(G)}.
 
 \mk
 
 b) Fix a Borel choice mapping ${\bf c}:\FF \to X$ \st for all $F\in \FF$, ${\bf c}(F)\in F$.
  Since for all $r>0$   the set 
$$\B=\{(F, a,b, r, J)\in \FF\times X^2\times {\mathbb Q}^+\times  \JJ: \  J\subset F,  \ e(J)=\{a,b\} \et   \diam(J)<r \}$$ 
is Borel  there exists a $\SSS$-measurable mapping ${\boldsymbol\tau}:\FF \times X^2 \times {\mathbb Q}^+\to \hat \JJ$ \st 
for all $(F,a,b,r)\in \FF\times X^2\times {\mathbb Q}^+$ if $J= {\boldsymbol\tau}(F,a,b,r)$ then:
$$\left\{
\begin{array}{ll}
(F,a,b,r, J)\in \B & \text{if} \  0<\d^F(a,b)<r\\
J=\{{\bf c}(F)\} & \text{if not}
 \end{array}\right.$$
So if we fix an enumeration $(m_p,n_p,r_p)_{p\in \wo}$ of the set $\wo^2\times {\mathbb Q}^+$ and
set for all $F\in \Carc^\Th(X)$:
 $${\boldsymbol\s_1}(F) (p)=  {\boldsymbol\tau} (F,  {\boldsymbol\s_0}(m_p),  {\boldsymbol\s_0}(n_p),r_p)$$
then the mapping  ${\boldsymbol\s_1}: \Carc^\Th(X) \to \hat \JJ^\wo$  satisfies condition 2) of Definition \ref{code}. It also follows from the hypothesis that the mapping  
${\boldsymbol\s}=({\boldsymbol\s_0},{\boldsymbol\s_1}):\Carc^\Th(X)  \to \A$ is  bianalytic  and for all   $F\in \Carc^\Th(X)$, $({\boldsymbol\s}(F), F)\in\C_\star$.
 Since ${\boldsymbol\s_0} $ is  $\G$-measurable and ${\boldsymbol\tau}$ is $\SSS$-measurable
 one easily checks that ${\boldsymbol\s_1}$, as well as  ${\boldsymbol\s}$, is  measurable relatively to the  $\s$-algebra class   generated by ${\mathcal A}(\G)$ and $\check{\mathcal A}(\G))$,
 hence ${\boldsymbol\s}$  is  ${\mathcal A}(\check{\mathcal A}( \G))$-measurable  and $\check{\mathcal A}({\mathcal A}( \G))$-measurable. It follows then   from  part a)  that  $\Carc^\Th(X)$ is in   $\check{\mathcal A}({\mathcal A}(\check{\mathcal A}(\G)))= {\bigl({\mathcal A}({\mathcal A}(\check{\mathcal A} }(\G)))\check{\bigr)}=\check{\mathcal A}(\check{\mathcal A}(\G))$. 
\end{proof}

  \begin{corol} \label{den} 
Let $X$ be a   Polish space.  If  $X$ contains at most countably many centers of triods then  the set $\Carc(X)$  is    in $\check{\mathcal A}(\check{\mathcal A}(\ca))$.      
\end{corol}

\begin{proof} 
Fix an enumeration $(c_n)_{n\in \wo}$ of the set $S^X$ of all centers of triods in $X$. 
Then by definition for all  $F\in\Carc^ \Th (X)$, the set $S^F$ is a non empty  subset  of $S^X$, hence
$\S^F\subset \adh{(\S^X\cap F)}^{\,\d^F}$.

Consider then the Borel  mapping
${\boldsymbol\s}_0: \Carc^ \Th (X)\to \wo^\wo$  defined by
$${\boldsymbol\s}_0(F)(k)=\left\{
\begin{array}{ll}
\min\{ n : \  c_n\in F  \et   \Vx j<k,  \ c_n \not = \s_0(F)(j),\} \} & \text{if this set is nonempty} \ \\
{\boldsymbol\s}_0(F)(k-1) &   \text{if not} 
\end{array} \right.$$
Note that if  $F\in \Carc^ \Th (X)$ then  $S^F\nv$  and so there exists $n$ \st   $c_n\in F$, so ${\boldsymbol\s}_0(F)(0)$ is  defined; and it follows by induction that 
${\boldsymbol\s}_0(F)(k)$  is  defined for all $k$. Consequently 
${\boldsymbol\s}_0$ is well defined and for all  $F\in \Carc^ \Th (X)$,  if ${\boldsymbol\s}_0(F)=\s$ then
$S^F\subset S\cap S^X= \{c_{\s (k)}; \ k\in \wo\}$,  hence  $\S^F\subset \adh{ \{c_{\s (k)}; \ k\in \wo\}}^{\d^F}$, and  ${\boldsymbol\s}_0$ satisfies the hypothesis of part b) of Theorem \ref{sG} with $\G=\ca$, 
so the set $\Carc(X)$  is    in $\check{\mathcal A}(\check{\mathcal A}(\ca))$.      
 \end{proof}

 \begin{corol} \label{sep} 
For any     Polish space  $X$, if  for all  $F\in \Carc^\Th(X)$   the set $\S^F$ is $\d^F$-separable,  then  the set $\Carc^\Th(X)$ is   ${\boldsymbol\Delta}^1_2$.    
\end{corol}

\begin{proof}  Note that  the set $\C^\star=\{(F,\a)\in  \FF\times\A: \ \SaF\supset \S^F\}$ is a $\ca$ subset of  
$\FF\times\A$,    
 and it follows from the hypothesis that $\pi(\C^\star)$ -- the projection of $\C^\star$ onto the first factor -- contains $\Carc^\Th(X)$.
By Kondo uniformization Theorem there exists  a  mapping  ${\boldsymbol\s}_0: \pi(\C^\star) \to \A$ with $\ca$ graph    contained in   $\C^\star$, hence satisfying
for all   $F\in \Carc^\Th(X)$, $ (F, {\boldsymbol\s}_0(F))\in  \C^\star$, so $ ({\boldsymbol\s}_0(F), F)\in  \C^\star$.
Since the graph of  ${\boldsymbol\s}_0$ is $\ca$ then   the mapping  ${\boldsymbol\s}_0$ is  
${\boldsymbol\S}^1_2$(hence ${\boldsymbol\Delta}^1_2$)-measurable on $ \Carc^\Th(X)$. Then 
 by Theorem \ref{sG} b), $\Carc^\Th (X)$  is    in $\check{\mathcal A}(\check{\mathcal A} ({\boldsymbol\Delta}^1_2)={\boldsymbol\Delta}^1_2$ by the invariance of   projective classes by operation $\mathcal  A$ (see \cite{kur} 38, Theorem 4).
 \end{proof}

\section{Back to planar Polish spaces}

  In all this section $X$ will be a {\em planar} Polish space. Note that  it follows from Corollary \ref{sep} and Theorem \ref{dsr} that the set 
  $\Carc(X)$ is   ${\boldsymbol\Delta}^1_2$, and our goal in this section is to prove the following more precise results.

\begin{thm} \label{codage} 
For any Polish space $X\subset \R^2$    there exists a $\SSS$-measurable mapping  {$\ss: \FF(X)\to \A$} \st for all $F\in \Carc^\Th( X)$,  $(\ss(F),F)\in \C_\star$.
\end{thm}

It follows then from Theorem \ref{sG}

\begin{corol}\label{upper} For any Polish space $X\subset \R^2$  the set $\Carc(X)$ is in  $\cospi$.
\end{corol}
  
 The proof of Theorem \ref{codage} relies strongly on the following notion.

\begin{ssect}  \label{P} Triods trap:  
{\rm  A {\em triod trap} is a quadruple $p=(W, I_0, I_1, I_2)$ where $W$ is the bounded connected component
of some Jordan curve  $\partial W$ and $I_0, I_1, I_2$ are three pairwise disjoint sub-arcs of $\partial W$; $W$ will be called the \textsl{domain} of $p$, and we set
 $ W=\dom(p)$  and $\diam (p)=\diam(W)$.

A triod $T$ will  said to be 
 {\em compatible} with the trap $(W, I_0, I_1, I_2)$ if    $T=J_0\cup J_1\cup J_2 \subset \adh W$   
 and for all $k$,   $J_k\cap \partial W=e(J_k)\mns \{{\bf c}(T)\}\subset I_k$.
 
 In fact  the heart of Moore's proof of Theorem \ref{moore}  (\cite{mo}, Lemma 1)  can 
 be restated as follows:

\begin{lem} \label{trap}{\sc (Moore)}. 
Two  triods    compatible with the same trap      have nonempty intersection.
 \end{lem}
To derive Theorem \ref{moore}  from   Lemma  \ref{trap}  observe that there exists a countable set $\P$ of  triod traps
 in the plane \st any triod in the plane is compatible with some $p\in \P$. 
For example one can take for $\P$   the countable set of  all \textsl{rational circular triod traps} i.e. 
the set of all traps of the form $p =(W, I_0, I_1, I_2)$  where $W$ is an open disc with rational radius and  rational coordinates center, 
and, for each $i\in3$, $I_i$  is  a sub-arc of $\partial  W$ with rational endpoints.

 \begin{figure}[ht ]
\begin{center}
\caption{}
\label{Trap}
\includegraphics[width=5cm]{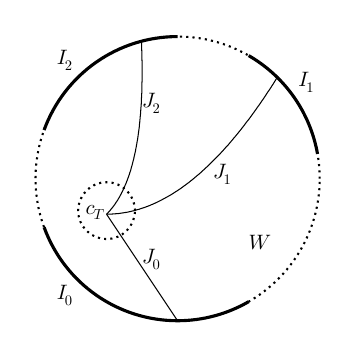}
\hskip 2cm
\includegraphics[width=5cm]{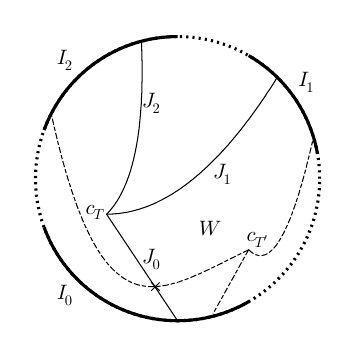}
\end{center}
\end{figure}

 Observe that if $T=J_0\cup J_1\cup J_2$ is compatible with the trap $p=(W, I_0, I_1, I_2)$, $W'$ is the bounded connected component
of some Jordan curve  $\partial W'$, $ \adh{ W'}\subset  W$ and $ {\bf c}(T)\in W'$ each $ J_k$ meets $ \partial W'$. It follows that  there is a unique  triod $ \tilde T=\td J_0\cup \td J_1\cup \td J_2$ contained in $T$ with $\td J_k\subset J_k$, hence  ${\bf c}(\td T) = {\bf c}(T)$, and   \st $\tilde J_k\cap \partial  W'$ is a singleton; and we shall
then write $\tilde T = T_{|W'}$. For fixed $ W$ and $ W'$  the set of such pairs $(T, \tilde T)$ is Borel, and the mapping $ T\mapsto \tilde T$ is Borel. 
If moreover $(I'_0, I'_1, I'_2)$ are pairwise disjoint sub-arcs of $\partial W'$ and $\tilde T$  is compatible with $p'=( W', I'_0, I'_1, I'_2)$ we shall  say that $ T$ is \textsl{weakly compatible} with the trap $p'$.
It follows that for all trap $p'$  the set of triods  $ T$ which are  weakly compatible with $ p'$ is Borel.

\sk
 }
 \end{ssect}

\nd{\it\underbar{Proof of Theorem \ref{codage}}:}
 
Let $\P$ denote the  set of  all  {rational circular triod traps} introduced in \ref{P}. For each triod trap $p $ the set $\T_p (X)$ of all triods in $X$ 
compatible with $ p $ is Borel, and the set $\M_p (X)$ of closed subsets
of $ X$ which match $p $, i.e which contain a triod compatible with $p $, is analytic.
For all $p\in \P$ and all $\eps>0$ we  denote by $ \P_{p,\eps}$ the set of all triod traps $ q\in \P$
 \st $ \adh{\dom( q)} \subset \dom (p)$ and $\diam(q) <\eps$.

\begin{lem}\label{join_centers}
For all $p\in \P $   there exists a Borel function $ \gamma_p $ which assigns to each pair of triods $(T, T')\in \TT_p (X)^2$   compatible with $p $ an arc joining $c_T$ to $ c_{T'}$ inside $ T\cup T'$.
\end{lem}
\begin{proof}
 By Lemma \ref{trap}, $ T\cap  T'\neq\0$. It is well known that there exists a Borel choice function $\rho$ assigning to each nonempty closed subset $F$ of a Polish space $ \X$ an element of  $ F$.
 So the function $ (T, T')\mapsto \rho( T\cap T')$ is Borel. Moreover there is a unique arc $ \vec J $ contained in $ T$ and joining $c_T$ to $ \rho(T\cap T')$ and similarly a unique arc $\vec J'$ 
 contained in $ T'$ and joining $c_{T'}$ to  $\rho( T\cap T')$. The graphs of these two mappings are Borel ; so the functions $ (T,T') \mapsto \vec J$ and $(T,T')\mapsto \vec J'$ are Borel, and so is the function 
 $(T,T')\mapsto  \vec J \vee \f(\vec J')$. And this latter arc joins $ c_T$ to $ c_{T'}$ inside $ T\cup T'$.
\end{proof}

\begin{lem}\label{phip}
 For each $p \in \P$ there exists a $\SSS$-meaurable function $ \ph_p $ defined on the set $\FF^*$ of all closed nonempty subsets of $X$ \st $\ph_p (F)=(T, a)$ 
 where $ T\in \TT_p (X)$, $T\subset F$  and $a = c_T$ 
 if $ F$ matches $p $ and $ \ph_p (F)=*$ if not.
\end{lem}
\begin{proof}
 The set $\tilde E_p =\{ (F, T, a)\in \FF^*\times \TT_p \times X : T\subset F \et a= c_T \et  F \text{ matches } T\,\}$ is Borel. Its projection $\M_p $ on $ \FF^*$ is analytic and 
 by Theorem \ref{YvN} there exists a $\SSS$-measurable selection of 
 $\tilde E_p $ on $\M_p $. Extending this  selection by $*$ on $ \FF^*\mns \M_p $ we define a $\SSS$-measurable function $ \ph_p $. For $ F\in \M_p $ we will set 
 $ \ph_p ^1(F)=  T$ and $ \ph_p ^2(F)= a$ if $ \ph_p (F) = (T,a)$.
\end{proof}

\begin{lem}
 Let $p\in \P$, $\eps>0$  and the $\SSS$-measurable function $\ph_p$ be as in Lemma \ref{phip}.  Then there exists a $\SSS$-measurable function $\lambda_{p,\eps} : \FF^*\to \P_{p,\eps}$
 \st for all $ F\in \FF^*$, if $ \ph_p(F) \neq *$,  the triod $\ph_p^1(F)$ is weakly compatible with the trap $ \lambda_{p,\eps}(F)$.
\end{lem}
\begin{proof}
 The set $\{F : \ph_p(F)\neq*\}$ is analytic hence belongs to $ \SSS$. 
 
 For each $ q\in \P_{p,\eps}$ the set $ \{ T\in \TT(X) : T \textsl{ is weakly compatible with }q\}$
 is Borel. Thus the set $\ZZ_q=\{ F : \ph_p(F) \textsl{ is weakly compatible with } q\}$ belongs to $\SSS$. And it is easy to check that $$
\M_p=\{F : \ph_p(F)\neq* \} \subset \bigcup_{q\in \P_{p, \eps}}\nolimits\ZZ_q$$
 since for all $ F\in \M_p$ the center of the triod $\ph^1_p(F)$ belongs to the domain $W'$ of some trap $ q=(W', I'_0, I'_1, I'_2)\in \P_{p, \eps}$.
 Thus it is possible to find a countable partition $ (\ZZ'_q)_{q\in \P_{p, \eps}}$ of $ \M_p$  \st $\ZZ'_q\in \SSS$ and $ \ZZ'_q\subset \ZZ_q$. Setting $$\lambda_{p,\eps}(F)=q \iff F\in \ZZ'_q$$
 completes the proof.
\end{proof}

\begin{lem}
 For all $(p , q )\in \P^2$, all $r\in \QQ^+$ and all $k\in \omega$ there exists a $\SSS$-measurable function $\psi_{p , q , r,k}$ on $ \FF^*$ \st 
 $\psi_{p , q , r,k}(F)$
  is an oriented arc $J$ joining $\ph^2_p(F)$ to $\ph^2_q(F)$ inside $F$
 with $\diam(J\cup \dom(p') \cup \dom(q') ) <r$ and $\max(\diam(p'), \diam(q'))<2^{-k}$ if such a  $J$ exists, and $\psi_{p , q , r,k}(F)=*$ if not.
\end{lem}
\begin{proof}
 As above the set
\begin{align*}
 \hat E_{p ,q , r}=\Bigl\{ (F, T, T', J)\in \FF^*\times \TT_p (X)&\times \TT_{q }(X)\times \JJ(X) :  J, T, T'\subset F\\
 &\et e(J)=(c_T, c_{T'})  \et \diam( J\cup p \cup q) <r \,\Bigr\}
 \end{align*}
 is Borel and has a $\SSS$-measurable selection $\chi_{p , q , r}$ on its projection $E^*_{p , q , r}=\pi( \hat E_{p ,q , r})$ on $ \FF^*$ 
 that we extend by $*$ on the coanalytic set 
 $\FF^*\mns E^*_{p , q , r}$.
 If $\chi_{p , q , r}(F)= (T, T', J)$ we denote $\chi_{p , q , r}^1(F)= T$, $\chi_{p , q , r}^2(F)= T'$ and $\chi_{p , q , r}^3(F)= J$.
 Then for $F\in E^*_{p , q , r}$, we necessarily have $ F\in \M_p  \cap \M_{q }$. 

 Clearly, if $F\in E^*_{p , q , r}$, we also have $F\in E^*_{p , q , r'}$ for all $ r'>r$ ; so the set 
 $\{(p , q , r) :F\in E^*_{p , q , r}\}$ is either empty or infinite.
 
 For $p,q\in \P$, $r\in \QQ^+$ and $k\in \omega$, we consider for all $ F\in \M_p\cap\M_q$ :  $p' = \lambda_{p, 2^{-k}}(F)$ , $ W'=\dom(p')$,  $q' = \lambda_{q, 2^{-k}}(F)$ and $ V'= \dom(q')$.
 
 By definition of $ \lambda_{p, 2^{-k}}(F)$ the triod $T=\ph_p^1(F)$ is weakly compatible with $p'$. So $T_{|W'}$ is compatible with $p'$, and
 ${\bf c}(T) ={\bf c}( T_{|W'})$. Similarly the triod $ S= \ph_q(F)$ is weakly compatible with $q'$.
 Thus
 using the  functions $\gamma_{p'}$ and $ \gamma_{q'}$ from Lemma \ref{join_centers}, we can consider the oriented arcs 
 \begin{align*}
 \vec I_{p, k,r, F}&= \gamma_{p'}\bigl(\ph^1_{p'} (F), \chi^1_{p', q' , r}(F)\bigr) \\
\vec J_{q, k, r,F}&= \gamma_{q'}\bigl(\ph^1_{q'} (F), \chi^2_{p', q' , r}(F)\bigr)
\end{align*}
respectively contained in $ W'=\dom(p')$ and $V'=\dom(q')$,  and joining respectively ${\bf c}(T) $ to $e_0( \chi^3_{p', q',r}(F))$ and 
${\bf c}(S) $ to $e_1( \chi^3_{p', q',r}(F))$.

Assume moreover that $F\in  E^*_{p,q,r}$.
Since  the functions $ \gamma_p$ are Borel, it is clear that the functions $ F\mapsto \vec I_{p,k,r,F}$ and  $F\mapsto \vec J_{q,k,r,F}$ are $ \SSS$-measurable,
and so is the concatenation $$
\psi_{p,q,r,k}(F)= \vec  I_{p,k, r,F}\vee \chi^3_{p', q' , r}(F) \vee \f(\vec J_{q, k, r,F})$$
which is an oriented arc joining $\ph^2_{p}(F)$ to $ \ph^2_{q}(F)$. Moreover $\psi_{p,q,r,k}(F) \subset \chi^3_{p', q' , r}(F) \cup W'\cup V'$, hence
$\diam(\psi_{p,q,r,k}(F)) <r$.
 \end{proof}
\mk

Applying Lemma \ref{compress} to the family $\bigl(\ph^2_p \bigr)_{p \in \P}$ and to the family $\bigl( \psi^3_{p , q , r}\bigr)_{(p , q , r)\in \TT(X)^2\times \QQ^+}$
we get $\SSS$-measurable functions $\tilde\ph : \FF^*\to X^\omega$ and $\tilde \psi : \FF^*\to \vec\JJ(X)^\omega$ such that for every $ F\in \FF^*$
\begin{align*}
\{ \tilde\ph_n(F) : n\in \omega\}& = \{ \ph^2_p (F) : F\in E_p , \ p \in \P\} \\
\{ \tilde\psi_n(F) : n\in \omega\} &= \{ \psi^3_{p , q , r,k}(F)   :  p , q \in \P, \ r\in \QQ^*,\ k\in \omega  \}
\end{align*}
Then the sequence $(a_n) = \bigl(\tilde\ph_n(F)\bigr)_n$ is $\tau$-dense in $ \Sigma^F$ since for every trap $p$ matched by $ F$
there is some $n$ \st $a_n$ is the center of some triod compatible with $p$. 

\begin{lem} \label{lem} 
If $(a_m, a_n)=(\tilde\ph_m(F),\tilde\ph_n(F)) \in E_F$ then
$$ \delta^F(a_m, a_n)= \inf\{ \diam (J) : \exists k\in \omega\  J=\tilde\psi_k(F) \et  e(J)=( a_m, a_n)\}\, .$$
\end{lem}

\begin{proof} 
 If  $(a_m, a_n)\in E_F$  then $\delta^F(a_n, a_m)<\infty$, and    for  all $r\in \QQ^+$ \st  $r> \delta^F(a_n, a_m)$  there exists an arc $J_0$ with endpoints $a_n$ and $ a_m$ \st $\delta^F(a_n, a_m)\leq \diam(J_0) <r$
 and an integer $k$ \st $\diam( J_0) + 2 ^{1-k} <r$. There exist $p$ and $ q$ in $\P$ \st $ a_n=\ph^2_p(F)$ and $ a_m=\ph^2_q(F)$, and we denote 
$p'= \lambda_{p, 2^{-k}}(F)$ and $ q'= \lambda_{p, 2^{-k}}(F)$. We then have $\diam(\vec I_{p,k,r,F}) < 2^{-k}$ and $\diam(\vec I_{q,k,r,F}) < 2^{-k}$, hence $$
\diam \bigl( \vec I_{p,k,r,F} \vee J_0 \vee \f(\vec I_{q,k,r,F}) \bigr) < 2^{-k} + \diam ( J_0) + 2^{1-k}  =\diam( J_0) + 2^{1-k} <r \, .$$
It follows that $ F\in E^*_{p',q', r}$ and that $ J_1=\psi_{p,q, r, F}$ satisfies $e(J_1)=( a_m, a_n)$ and $\diam( J_1) <r$. Then there exists $k\in \omega$ \st 
$J_1= \tilde\psi_k( F)$ and we are done.
\end{proof}

To finish the proof of Theorem \ref{codage}  define  $\ss: \FF(X)\to \C$   by 
$ \ss(F) = \bigl( \tilde \ph(F), \tilde \psi(F), F\bigr)$ and observe that 
if $F$ is arc-connected then any two elements $\tilde\ph_m(F),\tilde\ph_n(F)$ of $\tilde \ph(F)$ are $E_F$-equivalent  and apply Lemma \ref{lem}.
\hfill$\square$

\section{The exact complexity of $\Carc(\R^2)$}
\label{lowerBound}

By Corollary  \ref{upper} the set $\pcon$ of  arc-connected compact subsets of the plane $ \R^2$  
 is a $\cospi$ subset of $\KK( \R^2)$, and we now prove that this upper bound  complexity is optimal.

  \begin{thm}\label{lowerbound}
The set $\pcon$ is $\cospi$-complete.
\end{thm}

\begin{proof} 
We start by some preliminary constructions.
 
 \mk
 
\nd{\bf The Cantor space $\B$ :}\label{B}   
 We first define inductively for all $s\in \omega^{<\omega}$  reals $a_s$ and $ b_s$  in $\I=[0,1]$ as follows: 
 
 \sk

 Set  $a_\0=0$ ,  $b_\0=1$, and fix two increasing sequences $(a_n)_{n\in\wo}$ ,  
 $(b_n)_{n\in\wo}$ such that:
 $$a_\0 < a_0 < b_0 < a_1< b_1<\cdots a_n < b_n< a_{n+1} <\dots <b_\0 \ett b_\0= \sup _n a_n = \sup_n b_n$$
then for all $s\nv$ if 
$h_s:\xi\mapsto (1-\xi) a_s + \xi\, b_s$  is the   affine function such that 
 $h_s(\I)=
[a_s, b_s]$  define  $a_{s\cat \<n\>}= h_s(a_n)$ and $b_{s\cat \<n\>}= h_s(b_n)$.
\sk

Then  $ \rho=\sup_{n\in \omega} (b_n-a_n)<1$ and    for all  $s\in \mathbb S$, \ $b_s-a_s\leq \rho^{\abs s}$. Hence for all $ \sigma\in \omega^{\omega}$ 
   the real    $b_\s= \inf_{s\prec\sigma} b_s=\sup_{s\prec\sigma} a_s$,  is well defined, and we set:
$$ \B_0=\{b_s : s \in \omega^{< \omega}\} \quad ; \quad   \B_1=\{b_\s : \s \in \omega^{\omega}\}
\quad\text{and} \quad  \B=\B_0\cup \B_1$$
The  set  $\B$ is clearly a perfect compact subset  of  $\I$  with empty interior.

\sk

The next construction   is essentially due to Becker (see \cite{kc}, 33.17 and 37.11).
Let $\T\subset 2^{\seq}$ denote the set of all trees on $\omega$. 

\begin{lem}\label{becker}
There is a continuous function $B:\T \to \KK(\R^2)$ which assigns to any  $T\in\T$ a connected compact subset $B(T)$ of the unit square $\I^2$ \st :
\sn $a)$\quad   $(0,1)\in B(T) \et (\I\times\{0\})  \cup ( \{1\} \times \I ) \subset B(T) $,
\sn $b)$ \quad    $B(T)$ has at most two arc-components,
\sn    $\begin{aligned}
c) \quad  T \text{  is ill-founded }&\iff B(T) \text{ is arc-connected}\\
&\iff \text{there is an arc in $B(T)$ connecting } (0,1) \et (1,0). 
\end{aligned}$
\end{lem}

\begin{proof}

For any $u,v\in \R^2$ we  denote  by  $[u,v]$ the line segment joining $u$ to $v$. 
For all $s\in \seq$ let $\hat a_s,  \hat b_s$ the elements of $\R^2$ defined by  $\hat a_s=( a_s, 2^{-\abs s})$ and $\hat b_s=( b_s, 2^{1-\abs s})$, and consider the  family  $(R_s)_{s\in \seq}$ of compact subsets of $\I^2$ defined as follows:
$$
R_\0=( \I\times \{0\} ) \cup (\{1\}\times \I) \cup [\hat a_\0, \hat a_0] \cup \bigcup_{n\in \omega}\bigl( [\hat a_n, \hat b_n]\cup [\hat b_n, \hat a_{n+1}]\bigr)$$
 (see Fig.\;\ref{Beck}); 
 and  for all $s\in\seq$,   $ R_s$ is the image of $ R_\0$ under the affine mapping 
sending $\hat a_\0$ to $\hat a_s$ and $\{ b_\0\}\times \I$ to $\{ b_s\}\times [0, 2^{-\abs s}]$. 
Finally for all $ T\in \T$ let $ B(T)= \bigcup_{s\in T} R_s$.

\begin{figure}[ht] 
\begin{center}
\caption{}
\label{Beck}
\includegraphics[width=6cm]{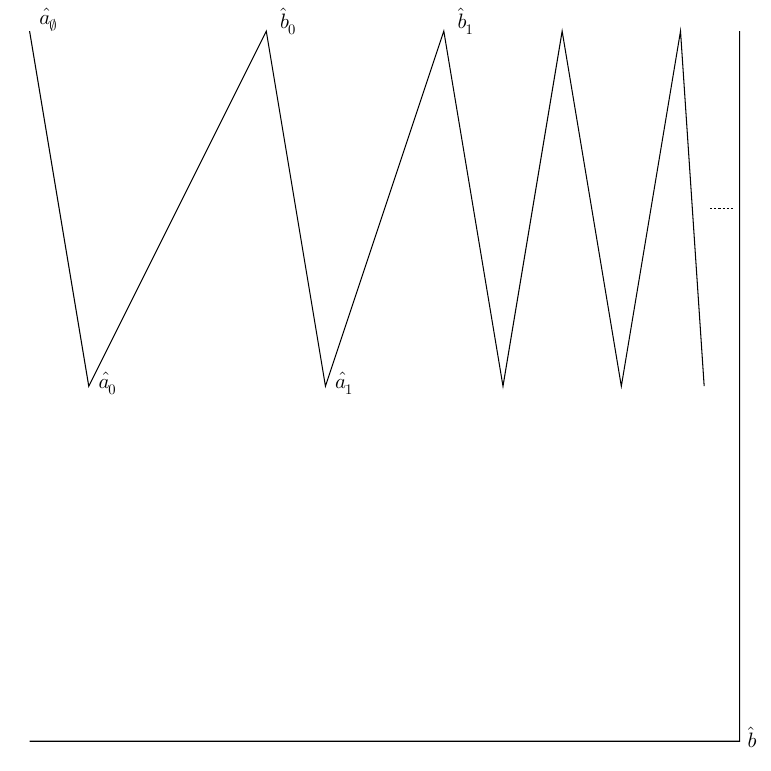}
\hskip2cm
\includegraphics[width=6cm]{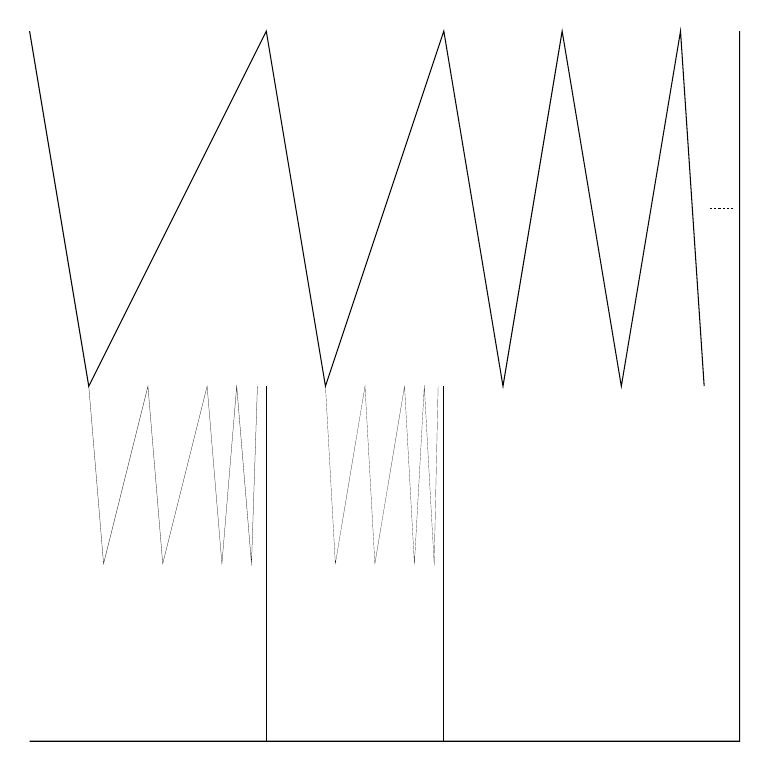}
\end{center}
\end{figure}

It  is clear that each $R_s$ is   a  connected, but not arc-connected, compact set. 
Since $ R_s$ and $R_{s\cat\<n\>}$ have $\hat a_{s\cat\<n\>}$ and $( b_s,0)$ in common, then  any  point in $B(T)$ is connected by an arc either to $\hat a_\0$ or to $\hat b=(1,0)$.
Hence $B(T)$ has at most two arc-components.

The set $B(T)$ is compact : indeed, if $(z_i)$ is a sequence in $ B(T)$ which converges to $ z\in \I^2$, there exists a sequence $ s^{(i)}\in T$ \st $z_i\in R_{s^{(i)}}$ 
and we can extract a subsequence \st 

\nd--~either $ s^{(i)}$ is a constant $s$, and then $ z\in R_s$ since $R_s$ is closed, 

\nd --~or there exists $s\in T$ and a sequence $ (n_i)$ 
converging to $\i$ \st $s\cat\< n_i\>\preceq s^{(i)}$, and then $z\in \{ b_s\}\times [0, 2^{-\abs s}]\subset R_s\subset B(T)$, 

\nd--~or else there exists   $\sigma\in \baire$ \st $\sigma_{|i}\preceq s^{(i)}$, 
and $z_i\to (b_\sigma, 0)\in R_\0\subset B(T)$.  

It is not difficult to see that $T\mapsto B(T)$ is continuous from $\T$ to $\KK(\R^2)$. Moreover if $\sigma$ is a branch of $ T$, then for all 
$s\prec\sigma$ there exists an arc $ J_n\subset R_s\subset B(T)$ with endpoints $\hat a_{\sigma_{|n}}$ and $\hat a_{\sigma_{|n+1}}$. And the concatenation of the $ J_n$'s yields 
an arc connecting $\hat a_\0$ to $(b_\sigma, 0)$. Thus in this case $B(T)$ is arc-connected. 

Conversely, if $B(T)$ is arc-connected let $ \gamma : \I\to B(T)$ be a continuous mapping with $\gamma(t)= \bigl(\gamma_1(t), \gamma_2(t)\bigr)$, 
$ \gamma(0)=\hat a_\0$ and $ \gamma(1)\in \I\times\{0\}$, and consider
  $\th=\inf\{ t : \gamma_2(t)=0\}>0$ and $\th_k=\sup\{ t<\th : \gamma_2(t) \geq 2^{-k}\}$. Then for $\th_k<t<\th$, $\gamma_2(t) < 2^{-k}$ 
and there exists some $s^{(k)}\in \omega^k$ \st 
$\gamma(t)\in \bigcup\{R_s; \ s\in T, \et s\succeq s^{(k)}\} $; so   $s^{(k)}\in T$ and  $ s^{(k)}\prec s^{(k+1)}$.
Hence   there exists $\sigma\in \baire$ \st $s^{(k)}\prec\sigma$ for all $k$. Thus $ \sigma$ is a branch of $ T$, and $ T$ is ill-founded.
\end{proof}

\nd\textbf{A connected compact set with uncountably many arc-components}

Let $\B$  the Cantor set constructed above with the elements   $a_s$, $b_s$ and $ b_\sigma$  for $s\in \seq$ and $\s\in \wo^\wo$, and fix    a decreasing sequence $(\alpha_n)\in \I$ \st $\alpha_0=1$ and  $\lim_n \alpha_n=\th =\dfrac23$.

\begin{figure}[ht] 
\begin{center}
\caption{}
\label{coBeck}
\includegraphics[width=6cm]{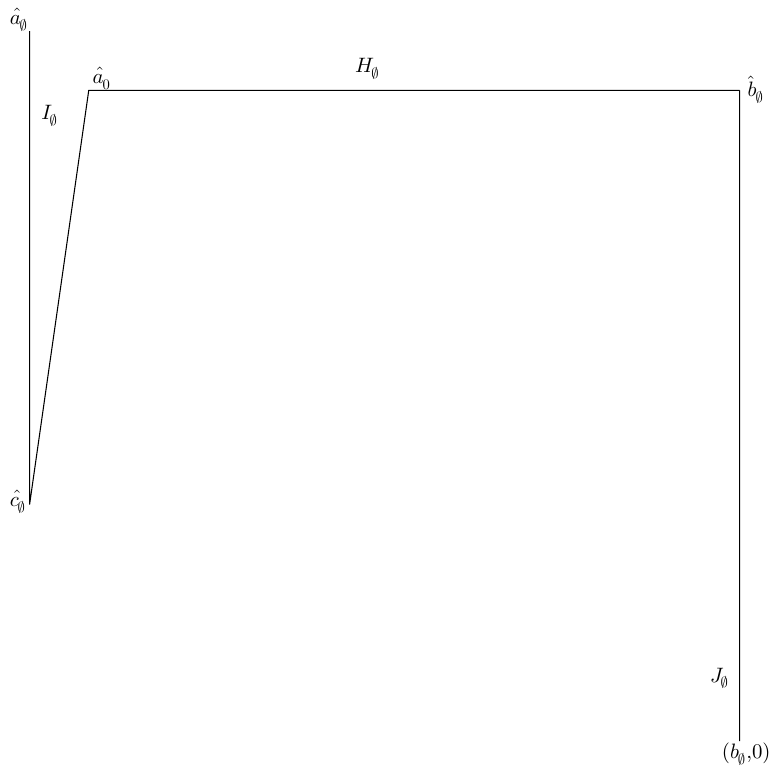}
\hskip2cm
\includegraphics[width=6cm]{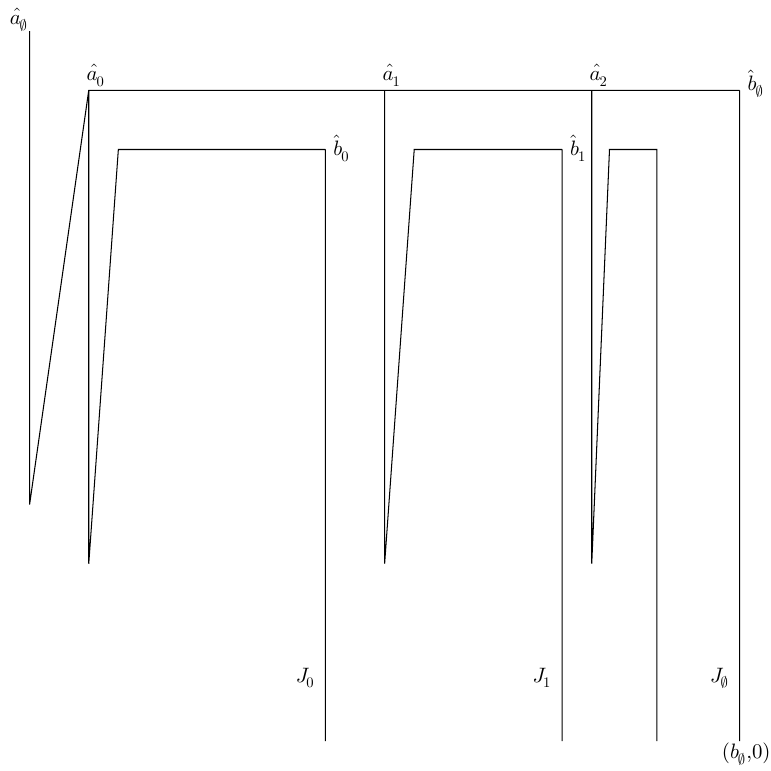}
\end{center}
\end{figure}

 For all $s\in \seq$  with $\abs s=k$,   
 let $P_s$ be the union of the following four segments of the unit square (see Fig. \ref{coBeck}): 
 
 -- the   vertical segment  $I_s$ with endpoints $\hat a_s= (a_s, \alpha_k)$ and 
$\hat c_s = (a_s,  \alpha_k -\th)$ (of length $\th$), 

-- the   segment  with  endpoints $\hat c_s$ and $\hat a_{s\cat\<0\>}$,

-- the horizontal segment $ H_s$ with endpoints 
$\hat a_{s\cat\<0\>}$ and $\hat b_s= ( b_s, \alpha_{k+1})$,

-- the  vertical segment $ J_s$ with endpoints $ \hat b_s$ and $(b_s, 0)$  
which has length $\alpha_{k+1}$. 

\nd It is clear that each $P_s$ is arc-connected  and has the point $\hat a_{s\cat\<n\>}$ in common with $ P_{s\cat\<n\>}$, hence $P_\i:= \bigcup_{s\in \seq} P_s$ is arc-connected, and
 the compact set $ P=\adh{ P_\i}$ is connected.
 
\begin{lem}\label{cobecker}
The set $P$ is the union of $P_\i$ and the vertical segments $ J_\sigma= \{ b_\sigma\}\times [0, \th]$. Moreover there is no arc in $P$ connecting $\hat a_\0$ to  any $ J_\sigma$.
\end{lem}

\begin{proof}
Let $(z_i)$ be a sequence of points of $ P_\i$ converging to $z=(x,y)\in P$. There exists a $ s^{(i)}\in \seq$ \st $z_i\in P_{s^{(i)}}$, and  again one can extract a subsequence \st:
 
-- either $ s^{(i)}$ is a constant $ s$ and then $z\in P_s\subset P_\i$ since $ P_s$ is closed, 

-- or there exists $s\in\seq$ and $(n_i)$ in $\omega$ converging to $\i$ \st $s\cat\<n_i\>\preceq s^{(i)}$, and then $z\in J_s\subset P_s\subset P_\i$, 
 
 -- or else there exists $\sigma\in \baire$ \st $\sigma_{|i} \preceq s^{(i)}$, and then $y\leq \th$ and $x=b_\sigma$, hence $z\in J_\sigma$. 
Conversely, if $z=(b_\sigma, y)\in J_\sigma$ , the points $ (a_s, y + \alpha_k-\th )$ for $ s=\sigma_{|k}$ belong to $P_\i$ and converge to $z$; hence  $J_\sigma\subset P$.

If there were an arc $J$ connecting $\hat a_\0$  and $ J_\sigma$, then $J$ should go through points in $ H_s$ for every $ s\prec\sigma$ 
and $J$ should contain $I_s$ for every $ s\prec\sigma$, which  is impossible. So each $ J_\sigma$ for $ \sigma\in \baire$ is an arc-component of $P$.
\end{proof}

\nd\textbf{Construction of a compact connected subset of $ \I^2$}

\nobreak
We recall that for any $s\in \seq$ with length $k=\abs s>0$ we denote by $ s^*$ the sequence of length $k-1$ \st $ s^*\prec s$.

We now modify the above constructed compact set $P$   by adding shortcuts between $ H_s$ and $ H_{s^*}$ :
choose for every $s\in \seq$ of length $k\geq 1$  some small square $Q_s$ inside the open rectangle $]a_{s\cat\<0\>}, b_s [ \times ]\alpha_{k+1}, \alpha_k[$ 
and consider the positive homothety $ h_s$ transforming the unit square $\I^2$ into $ Q_s$. For a given tree $ T_s\in \T$ the shortcut $ S_s$ will be 
the union of $h_s (B(T_s))$, where $ B(T_s)$ is the compact set defined in Lemma~\ref{becker},   with two vertical segments connecting respectively 
$ h_s(0,1)$ to $ H_{s^*}$ and $ h_s( 1,0)$ to $H_s$ (see Fig. \ref{xCoBeck}). 

\begin{figure}[ht]
\begin{center}
\caption{}
\label{xCoBeck}
\includegraphics[width=9cm]{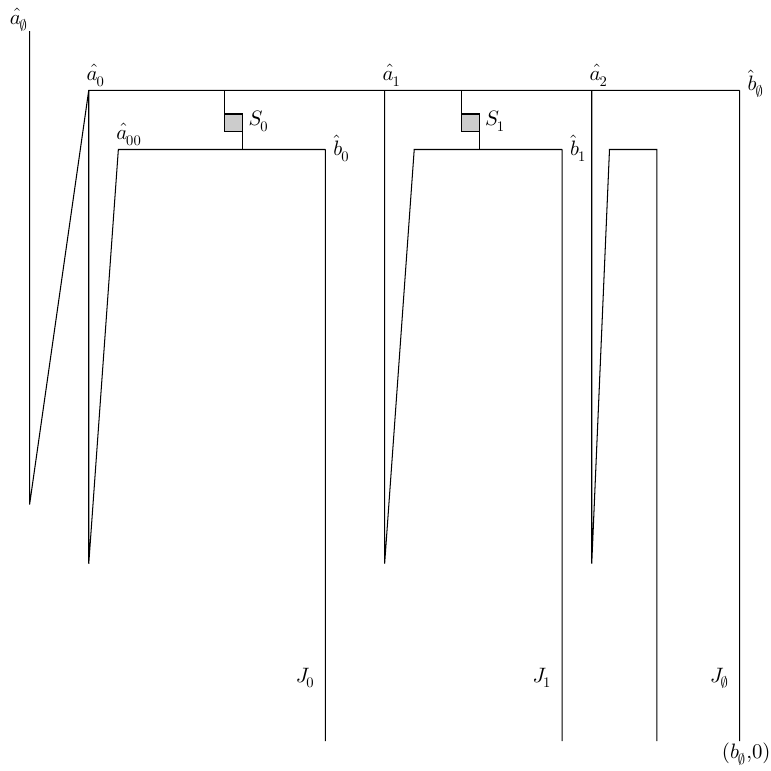}
\end{center}
\end{figure}

Then for any family 
$\tilde T=(T_s)_{s\in \seq}\in \T^{(\seq)}$, let $$
\Psi( \tilde T) = P\cup \bigcup_{k\geq 1} \nolimits \bigcup_{s\in\omega^k} \nolimits S_s.$$
\begin{lem}
  $\Psi(\tilde T)$ is compact and connected. Moreover the function $ \Psi$ is continuous from the compact space $ \T^{(\seq)}$ to $ \KK( \R^2)$.
\end{lem}
\begin{proof}
Since every point of $ B(T_s)$ is connected to $(0,1)$ or to $(1,0)$  (both if $ T$ is ill-founded), every point of $ S_s$ is connected to $ P_\i$. So $\Psi(\tilde T)$ 
is the union of the compact set $ P$ and the arc-connected set $P_\i\cup \bigcup_s S_s$, hence is connected.

To see that $\Psi(\tilde T)$ is compact, it is enough to prove that if a sequence $(z_i)$ in $ \bigcup_s S_s$ converges to $z=(x,y)\in \I^2$, then $ z\in \Psi(\tilde T)$. 
Then there exists a sequence $(s^{(i)}) $ in $\seq$ \st $ z_i\in S_{s_i}$. Again, if $s_i$ is constant equal to $ s$,  $S_s$ is closed and $z\in S_s\subset \Psi(\tilde T)$.
 If there exist $s\in \seq$ and $(n_i)$ tending to $\i$ \st $s\cat\<n_i\>\preceq s^{(i)}$, then $x=b_s$ and $ z\in J_s\subset P\subset \Psi(\tilde T)$. 
 And finally, if there exists $\sigma\in \baire$ \st $\sigma_{| i}\preceq s^{(i)}$, then $x= b_\sigma$ and $ y=\th$, hence $z\in J_\sigma\subset \Psi(\tilde T)$.

Observe   that for every $\eps>0$ there are only finitely many $s\in \seq$ \st $d(S_s, P) >\eps$.  It follows that  if we fix  an enumeration  $(s^{(n)})_n$  of $ \seq$, then
$ \Psi$ is the uniform limit of a sequence of continuous functions $\psi_n$ from $\T^{(\seq)}$ to $ \KK(\R^2)$ associating to any $ \tilde T$ the set $\psi_n(\tilde T)= P\cup \bigcup_{j=0}^n S_{s^{(j)}}$; hence  $\Psi$ is continuous.
\end{proof}
\mk

It results from what precedes that $\Psi(\tilde T)\mns \bigcup_{\sigma\in \baire} J_\sigma$ is arc-connected and that for each $\sigma\in \baire$, 
either $ J_\sigma$ is an arc-component of $ \Psi(\tilde T)$, or $J_\sigma$ is connected to $P_\i$ by an arc in $\Psi(\tilde T)$. 

\begin{lem}
For any $\sigma\in \baire$,  $J_\sigma$ is connected to $P_\i$ by an arc  in $\Psi(\tilde T)$ if and only if the set $ \{ s\prec\sigma : T_s\text{ is well-founded}\}$ is finite.
\end{lem}
\begin{proof}
Suppose first that there exists $ u\prec\sigma$ \st $T_s$ is ill-founded for each $s$ with $ u\prec s\prec \sigma$. Let $s^{(k)}\in\seq$ the beginning of
 $\sigma$ with  length $\abs u +k$.
Then since $ S_{s^{(k)}}$ is arc-connected  there exists for each $k\geq 1$ a continuous function 
$ \gamma_k : [1-2^{1-k}, 1-2^{-k}] \to  H_{s^{(k-1)}}\cup H_{s^{(k)}}\cup S_{s^{(k)}}$ \st 
$\gamma_k( 1-2^{1-k})= \hat a_{s^{(k-1)}\cat\<0\>}$ and $\gamma_k( 1-2^{-k})= \hat a_{s^{(k)}\cat\<0\>}$.  It follows that there exists a continuous function 
$\gamma :[0,1[ \to P_\i\cup \bigcup_s S_s$ which extends all $ \gamma_k$, that $ \gamma(0)\in P_\i$ and that $\gamma(\xi)\to (b_\sigma, \th)\in J_\sigma$ when $\xi\to 1$. 
This shows that $J_\sigma$ is connected to $P_\i$ by an arc in $\Psi(\tilde T)$.

Conversely, suppose that there exists a continuous path $\gamma :\I\to \Psi( \tilde T)$ connecting $\hat a_\0$ to $ J_\sigma$. 
Then $ \xi^*= \inf\{ \xi : \gamma(\xi)\in J_\sigma\}>0$.
It is easily seen that for any $s\in \seq$,  $ \Psi(\tilde T)\mns H_s$ is not connected and that $H_s$ separates $\hat a_s$ from $J_\sigma$ whenever $ s\prec\sigma$.
So $\xi_k= \sup\{ \xi\leq \xi^* : \gamma(\xi)\in H_{\sigma_{|k}}\}$ is well defined for all $k\geq 0$, and $\xi_k<\xi_{k+1}<\xi^*$. Then  if  $ \xi^{**}= \sup\xi_k \leq \xi^*$ 
we  necessarily have $\gamma(\xi^{**})\in \limsup_k H_{\sigma_{|k}}= \{ (b_\sigma, \th)\}\subset J_\sigma$, hence $\xi^{**}=\xi^*$. 
Moreover by continuity of $ \gamma$ at $\xi^*$ there exists $ k_0$ \st for $\xi_{k_0} <\xi<\xi^*$ : $d( \gamma(\xi), \gamma(\xi^*)) < 1-\th$.
 It follows that no point $\hat c_{\sigma_{|k}}$ for $ k>k_0$ can belong to $\gamma( [\xi_{k_0}, \xi^*])$. 
 Thus $\gamma$ has to go through all shortcuts $ S_{\sigma_{|k}}$ for $ k>k_0$. 
This implies that the corresponding  trees  $T_{\sigma_{|k}}$ are ill-founded for all $k>k_0$ and  the set $ \{ s\prec\sigma : T_s\text{ is well-founded}\}$ is finite.
\end{proof}
\begin{lem}
The compact set $ \Psi(\tilde T)$ is arc-connected if and only if for all $\sigma\in \baire$ there are only finitely many $T_s$ with ${s\prec\sigma}$  which are well-founded.
\end{lem}
\begin{proof}
The compact set $ \Psi(\tilde T)$ is arc-connected if and only if each $J_\sigma$ for $\sigma\in \baire$ is connected to $\hat a_\0$ by an arc;
and by the previous lemma this happens if and only if only finitely many $T_s$ with ${s\prec\sigma}$  are well-founded.
\end{proof}

\mk

To finish the  proof of Theorem \ref{lowerbound}  
let $Z\subset 2^\omega$ be any $\cospi$-set. We want to construct a continuous function $\Phi : 2^\omega\to \KK(\R^2)$ \st $\Phi(\zeta)$ is arc-connected if and only if $\zeta\in Z$. By definition of $\cospi$, there is a Souslin scheme $(\Gamma_s)_{s\in \seq}$ of $\ca$-sets, which can be assumed to be regular, \st 
$$
2^\omega\mns Z = \bigcup_{\sigma\in \baire} \bigcap_{s\prec\sigma} \Gamma_s
$$
Since the set $W\!F$ of well-founded trees is $\ca$-complete in $ \T$, we can find for each $s\in \seq$ a continuous function $T_s : \zeta\mapsto T_s(\zeta)$ \st 
$T_s(\zeta)\in W\!F\iff \zeta\in \Gamma_s$. Then the function $\tilde T : 2^\omega\to \T^{(\seq)}$ defined by $ \tilde T(\zeta) = (T_s(\zeta))_{s\in \seq}$ is  continuous too and so is 
$\Phi : \zeta\mapsto \Psi( \tilde T(\zeta))$.

If $\zeta\in Z$, then for all $\sigma\in \baire$, $\zeta\notin \bigcap_{s\prec\sigma} \Gamma_s$, thus there exists $s_0\prec\sigma$ \st $\zeta\notin \Gamma_{s_0}$ 
and since the Souslin scheme is regular we have also  for $ s_0\preceq s\prec\sigma$ : $\zeta\notin \Gamma_s$ and the tree $T_s(\zeta)$ is ill-founded;
 it follows that $ J_\sigma$ is connected to $\hat a_\0$ by an arc. Since this happens for all $ J_\sigma$, the compact set $ \Phi(\zeta)$ is arc-connected.

Conversely, if $\zeta\notin Z$ there is a $\sigma\in \baire$ \st for all $s\prec\sigma$ : $\zeta\in \Gamma_s$ and $ T_s(\zeta)\in W\!F$~;
 it follows that $J_\sigma$ is  connected to $\hat a_\0$ by none arc, and $ \Phi(\zeta)$ is not arc-connected.

 Thus $Z=\Phi\m(\pcon)$ and the proof is complete.
\end{proof}

  \bibliographystyle{amsplain}

  \end{document}